\documentclass[11pt]{article}
\usepackage{amsmath, amsthm, amssymb, amsfonts, amsbsy, color, graphics, epsfig, psfrag, graphicx, t1enc, mathtools, wasysym, enumitem}
\mathtoolsset{showonlyrefs}
\usepackage{appendix}
\usepackage[numbers,sort&compress]{natbib}
\usepackage[colorlinks=true]{hyperref}
\hypersetup{urlcolor=blue, citecolor=red}
\catcode`\@=11 \@addtoreset{equation}{section}

\catcode`\@=12

\newtheorem{theorem}{Theorem}[section]
\newtheorem{lemma}[theorem]{Lemma}

\newtheorem{proposition}[theorem]{Proposition}
\theoremstyle{definition}
\newtheorem{remark}{Remark}[section]

\DeclarePairedDelimiter{\abs}{\lvert}{\rvert}
\DeclarePairedDelimiter{\norm}{\lVert}{\rVert}

\DeclarePairedDelimiter{\skp}{\langle}{\rangle}
\DeclarePairedDelimiter{\prt}{(}{)}
\DeclarePairedDelimiter{\brk}{[}{]}

\newcommand{\ie}{\emph{i.e.}}
\newcommand{\cf}{\emph{cf.}\;}
\newcommand{\eg}{\emph{e.g.}}

\newcommand{\R}{{\mathbb R}} 

\newcommand{\PP}{{\mathcal{P}}}
\newcommand{\weakto}{\rightharpoonup}
\newcommand{\weakstar}{\stackrel{\ast}{\rightharpoonup}} 
\newcommand{\weak}{\mathrm{weak}}
\newcommand{\vrho}{\varrho}
\renewcommand{\phi}{\varphi}
\renewcommand{\rho}{\varrho}
\def\vec#1{\boldsymbol{#1}}
\newcommand{\vu}{\vec{u}}
\newcommand{\vU}{\vec{U}}
\newcommand{\vv}{\vec{v}}
\newcommand{\vf}{\vec{F}}
\newcommand{\vnu}{\vec{\nu}}

\newcommand{\M}{\mathcal{M}}
\newcommand{\E}{\mathcal{E}}
\newcommand{\KE}{\mathcal{KE}}
\newcommand{\divx}{{{\mathrm{div}}_x}\,}
\newcommand{\vr}		{\vrho}
\newcommand{\ep}		{\varepsilon}
\newcommand{\intO}{\int_\Omega}

\newcommand{\ol}{\overline}
\DeclareMathOperator{\Id}{\mathbb{I}_{d\times d}}
\DeclareMathAlphabet{\mathup}{OT1}{\familydefault}{m}{n}
\newcommand{\dx}[1]{\mathop{}\!\mathup{d} #1}
\newcommand{\ddt}{\frac{\dx{}}{\dx{t}}}
\newcommand{\mres}{\mathbin{\vrule height 1.6ex depth 0pt width
0.13ex\vrule height 0.13ex depth 0pt width 1.3ex}}

\title{Dissipative measure-valued solutions to the Euler-Poisson equation}
\author{Jos\'e A.\ Carrillo\footnotemark[1] \and Tomasz D\k{e}biec\footnotemark[2] \and Piotr Gwiazda\footnotemark[3] \and Agnieszka \'Swierczewska-Gwiazda\footnotemark[4]}

\begin{document}

\date{}
\maketitle

\abstract{We consider several pressureless variants of the compressible Euler equation driven by nonlocal repulsion-attraction and alignment forces with Poisson interaction. Under an energy admissibility criterion, we prove existence of global \emph{measure-valued solutions}, \ie, very weak solutions described by a classical Young measure together with appropriate concentration defects. We then investigate the evolution of a relative energy functional to compare a measure-valued solution to a regular solution emanating from the same initial datum. This leads to a (partial) weak-strong uniqueness principle.
}
\\[1em]
{\footnotesize
}

\renewcommand{\thefootnote}{\fnsymbol{footnote}}

\footnotetext[1]{Mathematical Institute, University of Oxford, Oxford, OX2 66G, UK.}

\footnotetext[2]{Laboratoire Jacques-Louis Lions, Sorbonne Universit\'e, 4 place Jussieu, 75005 Paris, France.}

\footnotetext[3]{Institute of Mathematics, Polish Academy of Sciences, \'Sniadeckich 8, 00-656 Warsaw, Poland.}

\footnotetext[4]{Institute of Applied Mathematics and Mechanics, University of Warsaw, Banacha 2, 02-097 Warsaw, Poland.}

\footnotetext[0]{Email addresses: carrillo@maths.ox.ac.uk, tomasz.debiec@sorbonne-universite.fr, pgwiazda@mimuw.edu.pl, aswiercz@mimuw.edu.pl.}



\section{Introduction}
Let $1\leq d\leq 3$, $\Omega\subset\R^d$ be a bounded smooth domain and $T>0$ be fixed.
We consider the following Euler-Poisson system with linear damping and confinement
	
\begin{equation}
	\begin{aligned}
	\partial_t \vr  + \divx (\vr \vu) & = 0,  &&\text{ in } (0,T) \times \Omega, 
	\\
	\partial_t (\vr \vu) + \divx (\vr \vu \otimes \vu) & =  - \gamma\vr \vu - \vr \nabla_x \Phi_\vr - \vr x,  &&\text{ in } (0,T) \times \Omega,          \label{eq:euler_1}
	\\
	-\Delta \Phi_\vr &= \vr - M_{\vr},  &&\text{ in } (0,T) \times \Omega,
	\end{aligned}
\end{equation}
where $\gamma>0$ is the friction coefficient and $M_\vr = \int_{\Omega}\rho(t,x)\dx{x}$ denotes the space average of the density. Note that, thanks to mass conservation, we have $M_\vr = \intO\rho_0(x)\dx{x}$, so that $M_\vr$ is constant in time.
The system is subject to the impermeability and Neumann boundary conditions

\begin{equation}\label{eq:bdryconditions}
    \vu\cdot n = 0,\qquad \nabla_x\Phi_\rho\cdot n = 0,\qquad\text{ on } (0,T)\times\partial\Omega,
\end{equation}
where $n$ denotes the outward unit normal to $\partial\Omega$.

\noindent To formulate the above equations in the measure-valued sense it is necessary to rewrite the nonlocal term $\vr \nabla_x \Phi_\vr$ into a divergence form. This can be done observing that any classical solution of~\eqref{eq:euler_1} satisfies the following pointwise identity 
\begin{equation}\label{eq:new_form}
	\prt*{\vr-M_\vr}\nabla_x\Phi_\vr = \frac{1}{2}\nabla_x  |\nabla_x \Phi_\vr |^2  - \divx[\nabla_x 	\Phi_\vr \otimes \nabla_x \Phi_\vr].
\end{equation}

\noindent Hence it is justified to consider the following form of the Euler-Poisson system instead of~\eqref{eq:euler_1}

\begin{equation}
	\begin{aligned}
	\partial_t \vr  + \divx (\vr \vu) & = 0, &\mbox{ in } (0,T) \times \Omega, 
	\\
	\partial_t (\vr \vu) + \divx (\vr \vu \otimes \vu) & =  - \gamma \vr \vu - \frac12\nabla_x  |\nabla_x \Phi_\rho |^2  + \divx[\nabla_x \Phi_\rho \otimes \nabla_x \Phi_\rho] - \vr x - M_\vr\nabla_x \Phi_\rho , &\mbox{ in } (0,T) \times \Omega,                               \label{eq:euler_2}
	\\
	-\Delta \Phi_\rho &= \vr - M_\vr, & \mbox{ in } (0,T) \times \Omega.
	\end{aligned}
\end{equation}

Multiplying formally the momentum equation of~\eqref{eq:euler_1} and using the continuity equation several times we arrive at the following identity
$$
    \partial_t\brk*{\frac12\vr|\vu|^2 + \vr\Phi_\vr + \frac12\vr|x|^2} - \divx\brk*{\prt*{\frac12\vr|\vu|^2 + \vr\Phi_\vr + \frac12\vr|x|^2 }\vu} = -\gamma\vr|\vu|^2 + \vr\partial_t\Phi_\vr,
$$
which, upon integration in space and using the Poisson equation, yields the following energy balance satisfied by smooth solutions of~\eqref{eq:euler_1}
\begin{equation*}
    \ddt\E(\vr,\vu,\Phi_\rho)(t) = -\gamma\intO\rho|\vu|^2\dx{x},
\end{equation*}
where
\[
    \E(\vr,\vu, \Phi_\rho) = \int_\Omega \frac{1}{2} \vr |\vu|^2 + \frac12|\nabla_x\Phi_\vr|^2 + \frac12\vr|x|^2 \dx{x},
\]
denotes the total energy associated with~\eqref{eq:euler_1}.

\subsection{Preliminaries and notation}

Our main results concern existence and conditional uniqueness of measure-valued solutions to the Euler-Poisson system~\eqref{eq:euler_2}, which will be defined in the following section. First however, we introduce the necessary notation and recall the formalism of Young measures.

\noindent Let $n,m\in\mathbb{N}$ and $X\subset\R^n, Y\subset\R^m$ be measurable sets. We denote the spaces of signed and non-negative Radon measures in $Y$ by $\M(Y)$ and $\M^+(Y)$, respectively. By $\PP(Y)$ we denote the space of probability measures. In each case we equip these spaces with the total variation norm
\begin{equation*}
	\norm{\nu}_{\mathrm{TV}} := \int_Y\dx{|\nu|}.
\end{equation*}

\noindent By $L^\infty_{\weak}(X;\M(Y))$ we denote the space of weakly-$^*$ measurable essentially bounded maps $\vnu=(\nu_x):X\to\M(Y)$. This means that for each $\phi\in C_0(Y)$ the map
\begin{equation*}
	x\mapsto \int_Y\phi(\lambda)\dx{\nu_x}(\lambda) \equiv \skp*{ \nu_x ; \phi}
\end{equation*}
is Lebesgue-measurable, where $C_0(Y)$ denotes the space of continuous real-valued functions on $Y$ which vanish at infinity. 
We note also that $L^\infty_{\weak}(X;\M(Y))$ is isometrically isomporhic to the dual space of the separable space $L^1(X;C_0(Y))$. This is a key point in associating a parameterised measure to a sequence of measurable functions. 

\medskip

When considering the issue of existence of solutions for system~\eqref{eq:euler_2}, we shall follow the usual strategy of first constructing weak solutions to an approximate problem, and then passing to the limit in the approximation parameters. Naturally, we run into the usual problem that the a priori estimates we can derive are not strong enough to allow for passing to the pointwise limit in the nonlinear terms. We must therefore find means to characterise the possible oscillation and concentration effects in our approximate sequences. This is done by embedding the problem into the larger space of bounded Radon measures:\ given a sequence $(z^\ep)_{\ep>0}$ of measurable functions $z^\ep:X\to Y$, we associate to each $z^\ep$ the mapping $\nu^\ep:X\to Y$ defined by $\nu^\ep(x) = \delta_{\{z^\ep(x)\}}$. Then the sequence $\prt*{\nu^\ep}_{\ep>0}$ belongs to the closed unit ball of $L^\infty_{\weak}(X;\M^+(Y))$. Therefore, by virtue of the Banach-Alaoglu Theorem, there exists a subsequence (which we shall never relabel) and a parameterised measure $\vnu = (\nu_x)_{x\in X}$ in $L^\infty_{\weak}(X;\M^+(Y))$ such that $\nu^{\ep}\weakstar\vnu$. In particular this implies that
	\begin{equation*}
		f(z^\ep) \weakstar \skp*{\nu_x ; f}
	\end{equation*}
in $L^\infty(X)$ for every $f\in C_0(Y)$. Moreover, clearly $\nu_x\geq 0$ and $\norm{\nu_x}_{\mathrm{TV}}\leq1$ for a.e.\ $x\in X$. The parameterised measure $\vnu\in L^\infty_{\weak}(X;\M^+(Y))$ is called the Young measure associated to (or generated by) the (sub)sequence $(z^\ep)$. The above observations, on various level of generality, are usually termed the Fundamental Theorem, see, \eg,~\cite{Ball1989, Pedregal}. The Young measure captures the oscillatory behaviour of a sequence and allows to characterise some nonlinear weak limits.

When working in the above setting and defining a measure-valued solution to a given problem, one usually desires that the Young measure describing the solution be a family of probability measures. It can be then thought of as giving the probability distribution of values of the physical quantities represented by the problem dependent variables around a given point in the physical space. To guarantee that this is the case some additional information on the underlying sequences is needed. More precisely, suppose that the sequence $z^\ep$ satisfies the following tightness condition:
\begin{equation}
\label{eq:tightness}
	\sup_{\ep>0}\int_X g(|z^\ep(x)|)\dx{x} < \infty,
\end{equation}
for some non-decreasing continuous function $g:[0,\infty)\to [0,\infty]$ with $\lim_{\alpha\to\infty}g(\alpha) = \infty$.
Then almost every of the $\nu_x$ is a probability measure. Furthermore, in this case one can show that for any $f\in C(Y)$ such that $(f(z^\ep))_{\ep>0}$ is weakly precompact in $L^1$, then
\begin{equation*}
	f(z^\ep) \weakto \skp*{\nu_x ; f}\quad \text{in $L^1(X)$}.
\end{equation*}
Notice that every sequence $(z^\ep)$, which is uniformly bounded in some Lebesgue space $L^p$, $1\leq p<\infty$, satisfies condition~\eqref{eq:tightness}.

Two important observations arise from the above discussion that are of importance in our approach. Firstly, when working with the variables $\rho, \vu, \nabla_x\Phi_\rho$, the estimates that we shall derive in Section~\ref{sec:existence} will depend on $\rho$. Therefore we will not be able to deduce any control over the velocity in the vacuum regions $\{\rho=0\}$. Thus, we shall be working with a parameterised family of positive measures, rather than probabilities. See also Remark~\ref{rem:othervariables}.
Secondly, due to lack of the pressure we cannot guarantee that the approximate quantities (in particular the density and the momentum) are uniformly integrable. Consequently, the nonlinearities that we deal with are neither $C_0$, nor are they weakly precompact in $L^1$. Therefore we have to introduce a way to capture possible concentration effects in the approximate sequences (which is done below), but also justify that the maps $x\mapsto\skp*{\nu_x;f}$ are integrable. This is formulated in Appendix~\ref{sec:appA}, together with an observation concerning projections of the Young measure onto individual variables.

Now let $f$ be a continuous function on $Y$, such that
\begin{equation*}
    \sup_{\ep>0}\int_X |f(z^\ep(x))|\dx{x} \leq C.
\end{equation*}
Then there is a measure $f_\infty\in\M(X)$, such that $f(z^\ep)$ converges (up to a subsequence) to $f_\infty$ weakly-$^*$ in $\M(X)$. One can then consider the difference
\begin{equation}
\label{eq:concentrationdefect}
    m^f = f_\infty - \skp*{ \nu_x ; f } 
\end{equation}
where $\nu_x$ is the Young measure generated by the sequence $(z^\ep)$.
We shall call this measure the \emph{concentration measure} associated to the function $F$. It is also sometimes called the concentration defect. Observe that $m^f = 0$ if the family $\prt*{f(z^\ep)}$ is weakly precompact in $L^1(X)$. As a note of caution let us point out that the term "concentration" might be misleading: the measure $m^f$ need not be supported on a set of small measure, see for instance~\cite{BallMurat1989} for an example of a sequence whose concentrations are being smeared out uniformly over the whole domain.

\noindent It is often needed to compare two concentration measures coming from two nonlinear functions. It turns out that if one of the nonlinearities dominates the other, then the same is true of the associated concentration measures. More precisely, we have the following result.

\begin{proposition}
\label{prop:concentrationrelations}
Let $\vnu=(\nu_x)$ be the Young measure generated by the sequence $(z^\ep)$.
If two continuous functions $f_1$ and $f_2\geq0$ satisfy $|f_1(z)|\leq f_2(z)$ for every $z\in Y$ and if $f_2(z^\ep)$ is uniformly bounded in $L^1(X)$, then we have
\begin{equation*}
    |m^{f_1}|(A) \leq m^{f_2}(A)
\end{equation*}
for every Borel set $A\subset X$.
\end{proposition}
\noindent The proof of this measure theoretic fact can be found in~\cite[Lemma~2.1]{FGSW2016} or~\cite[Proposition~3.3]{GKS2020}.

\medskip

In our case of interest we have $X = (0,T)\times\Omega$ and $Y=[0,\infty)\times\R^d\times\R^d$.
Since we will generally perform the same algebraic operations on both the oscillatory and concentration parts due to given nonlinearities of the Euler-Poisson equations, we will use the shorthand notation
\begin{equation*}
	\ol{f}(t,x) = \skp*{\nu_{t,x}(\lambda); f(\lambda)} + m^f(\dx{t},\dx{x}),
\end{equation*}
where $\lambda = (s,\vv,\vf) \in [0,\infty)\times\R^d\times\R^d$.

\noindent Let us make one final remark about concentration measures. As discussed above, we generally have $m^f\in\M(X)$. It is also not difficult to show that if $f$ is non-negative, then so is the corresponding concentration measure. 
When $X = (0,T)\times\Omega$, then oftentimes it is also possible to guarantee a disintegration of $m^f$ with respect to the time and space variables. 
In Section~\ref{sec:existence} we prove existence of dissipative measure-valued solutions using a sequence of weak solutions to an approximate problem. Due to the energy inequality we will obtain bounds in $L^\infty(0,T;L^1(\Omega))$ for the approximate quantities. Then the corresponding concentration measures will admit a disintegration of the form
\begin{equation*}
	m^f = \dx{m}^f_t(x)\otimes \dx{t}, 
\end{equation*}
where the family $t\mapsto m^f_t$ is bounded and weakly-$^*$ measurable.

\subsection{Main results and structure of the paper}

We state below our main result concering the Euler-Poisson system~\eqref{eq:euler_2}. It is a weak-strong type result comparing a measure-valued solution and a regular solution emanating from the same finite-energy initial data. By a "regular" or "strong" solution we shall mean a continuously differentiable triple $(\rho,\vU,\Phi_r)$ with bounded velocity, which satisfies either (and thus both) system~\eqref{eq:euler_1} or~\eqref{eq:euler_2} pointwise. The measure-valued solutions will be defined only in the following section, however for now let us just say that we they comprise of a classical Young measure $\vnu$ together with a number of concentration-defect measures as described above, which satisfiy equations~\eqref{eq:euler_2} in an averaged sense. Furthermore, they are required to be "dissipative" or "admissible" by exhibiting an energy inequality, see Section~\ref{sec:dmvs}. Vitally, the following theorem concerns precisely the class of measure-valued solutions for which we can show global existence.

\begin{theorem}\label{thm:mv-stronguniqness}
 Let $1\leq d\leq 3$ and $\Omega\subset\R^d$ be a bounded smooth domain. Let 
    \[
        (r,\vU,\Phi_r)\in C^1([0,T)\times\bar{\Omega};(0,\infty))\times C^1([0,T)\times\bar{\Omega};\R^d)\times C^2([0,T)\times\bar{\Omega})
    \]
    be a regular solution of~\eqref{eq:euler_2} with initial data $r(0,x)=r_0(x),\, \vU(0,x)=\vU_0(x)$ of finite energy and let $(\vnu, m^\rho, m^{\rho\vu}, m^{\rho \vu\otimes \vu}, m^{|\nabla\Phi|^2}, m^{\nabla\Phi\otimes\nabla\Phi})$ be a dissipative measure-valued solution to the system~\eqref{eq:euler_2} with initial state
\[
    \nu_{0,x} = \delta_{\{r(0,x),\vU(0,x),\nabla\Phi_r(0,x)\}},\;\;\;\; \text{for a.e.}\;\; x\in\Omega.
\]
Then 
\[
m^{\nabla\Phi\otimes\nabla\Phi}=0,\;\; m^{|\nabla\Phi|^2}=0,
\]
and we have the following identifications
\begin{align*}
	\skp*{\nu_{t,x} ; \rho} + m^{\rho} &= r,\\	
	\skp*{\nu_{t,x} ; \rho\vu} + m^{\rho\vu} &= r\vU,\\
	\skp*{\nu_{t,x} ; \rho\vu\otimes\vu} + m^{\rho\vu\otimes\vu} &= r\vU\otimes\vU,\\
	\skp*{\nu_{t,x} ; \rho|\vu|^2} + m^{\rho|\vu|^2} &= r|\vU|^2,
\end{align*}
which hold for almost every $(t,x)\in(0,T)\times\R^d$.
Furthermore, the Young measure admits the decomposition
\begin{equation*}
	\nu_{t,x} = \bar{\nu}_{t,x} \otimes \delta_{\{\nabla\Phi_r(t,x)\}},
\end{equation*}
for some parameterised measure $\bar{\nu}\in L^\infty_{\weak}((0,T)\times\R^d;\M^+([0,\infty)\times\R^d))$; and in turn the restriction $\bar{\nu}\mres ((0,\infty)\times\R^d)$ decomposes into
\begin{equation*}
	\bar{\nu}_{t,x} \mres ((0,\infty)\times\R^d) = \bar{\bar{\nu}}_{t,x} \otimes \delta_{\{\vU(t,x)\}}
\end{equation*}
for some parameterised measure $\bar{\bar{\nu}}\in L^\infty_{\weak}((0,T)\times\R^d;\M^+(0,\infty))$.
Finally, all the concentration measures $m^{\rho}$, $m^{\rho\vu}$, $m^{\rho\vu\otimes\vu}$ and $m^{\rho|\vu|^2}$ are absolutely continuous with respect to the Lebesgue measure.
\end{theorem}

The above theorem is not a usual weak-strong uniqueness result, where we would like to assert that the Young measure is necessarily a tensor product of Dirac masses concentrated at the strong solution and all the concentration measures vanish.
Indeed, this was the case for most of the recent studies on weak-strong uniqueness for measure-valued solutions in hydrodynamics or more general conservation laws. Starting from the works of Brenier et al.~\cite{BrDeLeSz2011} on the incompressible Euler and Demoulini et al.~\cite{Demoulini2012} on polyconvex elastodynamics, and their somewhat surprising observation that measure-valued solutions can enjoy the weak-strong uniqueness property (under an admissibility condition), this property has been proved for a variety of other equations, see~\cite{GSW2015, GKS2020, FGSW2016, ChristoforouTzavaras2018}. Notably the weak-strong uniqueness property for dissipative measure-valued solutions has been recently put to practical use in proving convergence of finite volume numerical schemes for the Euler and Navier-Stokes equations~\cite{Feireisl:2018sy, Feireisl:2020fj}.
The common feature of these results is that one can control each relevant term in the equations by a relative energy between the measure-valued and regular solution in such a way that the vanishing of the relative energy implies that the corresponding measure-valued and regular quantities have to coincide.
In the situation considered in this paper we cannot quite follow the same algorithm:\ establishing a relative energy inequality does not immediately imply a weak-strong uniqueness result.
This is due to lack of appropriately strong information on the density and is a fundamental feature of system~\eqref{eq:euler_2} due to its pressurelessness. One might like to compare this situation with the isentropic compressible Euler equations with pressure $p(\rho)\sim\rho^\gamma$, see~\cite{GSW2015}. The potential energy term of the corresponding relative energy functional:
\[
    E_{rel}^{mv} = \intO \frac12\ol{\rho|\vu-\vU|^2} + \ol{\frac{1}{\gamma-1}\rho^\gamma-\frac{\gamma}{\gamma-1}r^{\gamma-1}\rho+r^\gamma}\dx{x}
\]
being zero immediately implies (due to convexity of the pressure) that the projection of the Young measure solution onto the first variable equals $\delta_{\{r(t,x)\}}$. This information can be used in the first term of the relative energy to conclude that the measrue-valued solution decomposes into $\nu_{t,x} = \delta_{\{r(t,x)\}}\otimes\delta_{\{\sqrt{r(t,x)}\vU(t,x)\}}$ almost everywhere. This reasoning cannot be mimicked in the pressureless case.

Let us further mention that the above issue is also related to some classical results on compensated compactness in the one-dimensional isentropic gas dynamics, see for instance~\cite{LPT94a, DiPerna1985, LPS1996}. Given $L^\infty$ initial data $\rho_0, u_0$ and a sequence $\rho_n, u_n$ of entropy weak solutions to the one-dimensional Euler equations, one can show that there holds the convergences
\begin{equation*}
	\rho_n \rightarrow \rho,\quad u_n\rightarrow u
\end{equation*}
to another entropy solution in the almost everywhere sense, where the second convergence holds only on the set $\{\rho(t,x)>0\}$. In terms of Young measures this means that on this set the Young measure reduces to Dirac masses, \ie,
\begin{equation*}
	\nu_{t,x} = \delta_{\{\rho(t,x)\}}\otimes \delta_{\{u(t,x)\}},\quad\text{if }\;\rho(t,x)>0;
\end{equation*}
while on the vacuum set $\{\rho(t,x)=0\}$ one has
\begin{equation*}
	\nu_{t,x} = \delta_{\{\rho(t,x)\}} \otimes \bar{\nu}_{t,x},
\end{equation*}
for some probability measure $\bar{\nu}$ whose support is contained in a compact set $[-a,a]$ depending on the initial data through the so-called Riemann invariants and invariant sets (in the $L^\infty$ sense). Thus, in this case the invariant sets prevent concentrations and guarantee the tightness condition in the "$u$-direction" (which is not necessarily true in general). Notice that the use of Riemann invariants is only available in one spatial dimension.

The paper is structured as follows. We first introduce the notion of dissipative measured-valued solutions to the Euler-Poisson system in the next section. We analyse the existence of these solutions for the Euler-Poisson system in section 3. Section 4 deals with the relative entropy argument between strong solutions and any dissipative measured-valued solution. We finally prove our partial result about weak-strong uniqueness of these solutions in Section 5. Section 6 is finally devoted to a generalization to the case of adding alignment terms to the Euler-Poisson system as in the case of Cucker-Smale models. 


\section{Dissipative measure-valued solutions to the Euler-Poisson equation}\label{sec:dmvs}

In what follows we use dummy variables $(s, \vv, \vf)\in[0,\infty)\times\R^d\times\R^d$ when integrating with respect to a parameterised measure $\vnu\in L^\infty_{\weak}( (0,T) \times \Omega; \M^+([0,\infty) \times \R^d \times \R^d))$. They should be thought of as representing $\rho$, $\vu$ and $\nabla_x\Phi_\rho$, respectively.

We say that $(\vnu, m^\rho, m^{\rho\vu}, m^{\rho \vu\otimes \vu}, m^{|\nabla\Phi|^2}, m^{\nabla\Phi\otimes\nabla\Phi})$ is a dissipative measure-valued solution of the Euler-Poisson system \eqref{eq:euler_2} in $(0,T) \times \Omega$ with initial data $(\vnu_0, m^\rho_0, m^{\rho\vu}_0, m^{\rho \vu\otimes \vu}_0, m^{|\nabla\Phi|^2}_0, m^{\nabla\Phi\otimes\nabla\Phi}_0)$ if
\[
    \vnu = \prt*{\nu_{t,x}}_{(t,x) \in (0,T)\times \R^d} \in L^\infty_{\weak}\prt*{(0,T) \times \Omega; \M^+ \prt*{[0,\infty) \times \R^d \times \R^d}}\;\;\text{is a parameterized measure},
\]
\[
    m^\rho, m^{|\nabla\Phi|^2} \in L^\infty(0,T;\M^+(\overline\Omega)),\;\; m^{\rho\vu} \in  L^\infty(0,T;\M(\overline\Omega)^d),\;\; m^{\rho \vu\otimes \vu}, m^{\nabla\Phi \otimes \nabla\Phi} \in L^\infty(0,T;\M(\overline{\Omega})^{d\times d}),
\]
and the following hold:

\begin{itemize}[label={}]

\item 
{\emph{Continuity equation}}:
	\begin{equation}\label{eq:mvcontinuity}
	\begin{split}
	\int_{\Omega}  & \skp*{ \nu_{\tau,x} ; s } \psi (\tau, x ) \dx{x} - 
	\int_{\Omega} \skp*{ \nu_{0,x}; s } \psi(0, x) \dx{x} + \intO\psi(\tau,x)\dx{m}^{\rho}_\tau(x) - \intO\psi(0,x)\dx{m}^{\rho}_0(x)
	\\[0.5em]  
	& = \int_0^\tau\!\! \int_{\Omega} 
	\brk*{
	\skp*{ \nu_{t,x} ; s } \partial_t \psi + \skp*{ \nu_{t,x} ; s \vv } \cdot \nabla_x \psi} \dx{x}\dx{t}
	+ \int_0^\tau\!\!\intO\partial_t \psi\dx{m}^\rho(x)\dx{t} + \int_0^\tau\!\!\intO\nabla_x\psi\cdot\dx{m}^{\rho\vu}(x)\dx{t},
	\end{split}
	\end{equation}
for a.e.\ $\tau \in (0,T)$ and every $\psi \in C^1 ([0,T] \times \overline{\Omega})$;

\item
{\emph{Momentum equation}}:
	\begin{equation}\label{eq:mvmomentum}
	\begin{split}
	    \int_{\Omega} & \skp*{ \nu_{\tau,x} ;  s \vv } \cdot \phi (\tau, x) \dx{x} 
	    - \int_{\Omega} \skp*{ \nu_{0,x}; s \vv } \cdot \phi (0, x) \dx{x} 
	    + \intO\phi(\tau,x)\cdot\dx{m}^{\rho\vu}_\tau(x) 
	    - \intO\phi(0,x)\cdot\dx{m}^{\rho\vu}_0(x)
	    \\[0.5em] 
	    & =  \int_0^\tau\!\! \int_\Omega 
	    \skp*{ \nu_{t,x}; s \vv }  \cdot \partial_t \phi \dx{x}\dx{t}
	    + \int_0^\tau\!\! \int_\Omega  \skp*{ \nu_{t,x}; s \vv\otimes\vv } : \nabla_x \phi \dx{x}\dx{t}
	    - \gamma\int_0^\tau\!\! \int_\Omega\skp*{ \nu_{t,x}; s \vv } \cdot \phi\dx{x}\dx{t} 
	    \\[0.5em]
	    &\hspace{0.2cm}
	    - \int_0^\tau\!\! \int_\Omega\skp*{ \nu_{t,x}; s } x\cdot\phi\dx{x}\dx{t}
	    + \int_0^\tau\!\! \int_\Omega\frac12\skp*{\nu_{t,x}; |\vf|^2}\divx\phi \dx{x}\dx{t} 
	    - \int_0^\tau\!\!\int_\Omega \skp*{ \nu_{t,x};\vf\otimes\vf}:\nabla_x\phi \dx{x}\dx{t} 
	    \\[0.5em]
	    &\hspace{0.2cm}
	    - M_\rho\int_0^\tau\!\!\intO\skp*{\nu_{t,x};\vf}\cdot\phi \dx{x} \dx{t}
	    +\int_0^\tau\!\!\intO\partial_t \phi \cdot \dx{m}^{\rho\vu}(x)\dx{t} 
	    + \int_0^\tau\!\!\intO\nabla_x \phi : \dx{m}^{\rho \vu\otimes \vu}(x)\dx{t} 
	    \\[0.5em]
	    &\hspace{0.2cm} 
	    - \gamma\int_0^\tau\!\!\intO\phi \cdot \dx{m}^{\rho\vu}(x)\dx{t}
	    -\int_0^\tau\!\!\intO \phi\cdot x \dx{m}^{\rho}(x)\dx{t}
	    +\int_0^\tau\!\!\intO\divx\phi\dx{m}^{\abs*{\nabla\Phi}^2}(x)\dx{t} \\[0.5em]
	   &\hspace{0.2cm}
	    -\int_0^\tau\!\!\intO\nabla\phi : \dx{m}^{\nabla\Phi\otimes\nabla\Phi}(x)\dx{t},
	\end{split}
	\end{equation}
for a.e.\ $\tau \in (0,T) $ and every $\phi \in C^1([0,T] \times \overline{\Omega}; \R^d)$;

\item 
\emph{Poisson equation}: 
For a.ae.\ $\tau \in (0,T)$ and every $\vartheta \in C^1 (\overline{\Omega})$
	\begin{equation}\label{eq:mvPoisson}
		\int_\Omega \skp*{ \nu_{\tau,x}; \vf}\cdot\nabla_x\vartheta\dx{x} = \int_\Omega\skp*{\nu_{\tau,x}; s}\vartheta \dx{x} + \intO\vartheta(x)\dx{m}^{\rho}_\tau(x) - M_\vr\int_\Omega\vartheta(x)\dx{x};
	\end{equation}
and furthermore there exists a function $\Psi_\rho\in L^2(0,T; W^{1,2}(\Omega))$ such that	
	\begin{equation}\label{eq:gradientmeasure}
    	\nabla_x\Psi_\rho(t,x) = \intO \vf \dx{\nu_{t,x}}(s,\vv,\vf),
	\end{equation}
almost everywhere.

\item 
{\emph{Energy inequality}}:
There is a non-negative measure $m^{\rho|\vu|^2}\in L^\infty(0,T;\M^+(\overline{\Omega}))$ such that
	\begin{equation}\label{eq:energy_inequality}
	\begin{split}
	\int_\Omega& \skp*{\nu_{\tau,x}; \frac{1}{2} s | \vv|^2} \dx{x} + \frac12 m_\tau^{\rho|\vu|^2}(\Omega) + \int_\Omega\skp*{ \nu_{\tau,x}; \frac12|\vf|^2} \dx{x} + \frac12m_\tau^{|\nabla\Phi|^2}(\Omega) + \int_\Omega\skp*{ \nu_{\tau,x}; s} \frac12|x|^2 \dx{x}
	\\[0.5em]
	&\qquad+\intO\frac12|x|^2\dx{m}^{\rho}_\tau(x)
	\\[0.5em]
	&\leq \int_\Omega\skp*{ \nu_{0,x}; \frac{1}{2} s | \vv|^2 } \dx{x} + \frac12 m_0^{\rho|\vu|^2}(\Omega) + \int_\Omega\skp*{ \nu_{0,x}; \frac12|\vf|^2} \dx{x} + \frac12m_0^{|\nabla\Phi|^2}(\Omega) + \int_\Omega\skp*{ \nu_{0,x}; s} \frac12|x|^2 \dx{x} 
	\\[0.5em]
	&\qquad+\intO\frac12|x|^2\dx{m}^{\rho}_0(x) - \int_0^\tau \brk*{\int_{\Omega} 
	\skp*{ \nu_{t,x}; s |\vv|^2 } \dx{x} + m^{\rho|\vu|^2}(\Omega)}\dx{t},
	\end{split}
	\end{equation}
for a.e.\ $\tau \in (0,T)$.
\end{itemize}

\begin{remark}\hfill
    \begin{enumerate}

		\item We shall usually think of a measure-valued solution as being a weak limit of some family of weak solutions (either to an approximate problem or to~\eqref{eq:euler_2} itself). Then all the concentration-defect measures are easily described by~\eqref{eq:concentrationdefect} and can be related to the formalism of generalised Young measures via so-called recession functions, see~\cite{AlBo}. Moreover Proposition~\ref{prop:concentrationrelations} gives natural relationships between these measures, which are inherited from the corresponding inequalities at the approximate level. In particular all the concentration-defects appearing in the governing equations are dominated by those appearing in the energy inequality. Notice however that there might exist measure-valued solutions which do not arise as limits of weak solutions, see~\cite{Chiodaroli:2017xy, GallenmullerWiedemann2021}.
    
        \item Let us note that the appearance of a non-trivial concentration measure in the density cannot be excluded, because of low integrability of $\rho$. This is in contrast to the situation, when a pressure term is present in the system. Then typically $p(\rho)\sim \rho^\gamma$, $\gamma>1$, implying uniform intergrability, and excluding the possibility of concentrations in the density.
        
        \item Notice that one cannot exclude concentrations in the term $\intO\frac12|\nabla_x\Phi_\rho|^2\dx{x}$ in the energy, since the Poisson equation $-\Delta\Phi_\rho=\rho\in\mathcal{M}_{+}$ only guarantees $\nabla_x\Phi\in L^p$ with $\frac{d}{d-1}<p\leq2$.
        
        \item Let us recall the shorthand notation
        $$
            \ol{f}(t,x) = \skp*{ \nu_{t,x}(\lambda); f(\lambda) } + m^f(\dx{t},\dx{x}).
        $$
        Thus, equation~\eqref{eq:mvcontinuity} can be written concisely as
        \begin{equation*}
	        \int_{\Omega} \ol{\rho}(\tau,x)\psi(\tau,x)\dx{x} - 
	        \int_{\Omega} \ol{\rho}(0,x)\psi(0, x) \dx{x}
	        = \int_0^\tau \int_{\Omega} 
	        \brk*{
	        \ol{\rho}\partial_t \psi + \ol{\rho\vu}\cdot\nabla_x \psi} \dx{x}\dx{t},
        \end{equation*}
        and similarly for the other equations.
        
        \item The constant $M_\vr$ in equations~\eqref{eq:mvmomentum} and~\eqref{eq:mvPoisson} should more precisely be written as $M_{\bar{\vr}}$ to account for possible concentrations in the initial data, when one has
        \[
        M_{\bar{\vr}} = \intO\skp*{\nu_0;s}\dx{x} + m_0^\vr(\Omega).
        \]
        However, in the proof of existence, as well as the subsequent analysis, we will specialise to the case of initial data without concentrations. 
        Let us also observe that testing the continuity equation~\eqref{eq:mvcontinuity} with the test function $\psi\equiv 1$, we obtain
        \[
        \intO \skp*{\nu_{\tau,x};s}\dx{x} + m_\tau^{\vr}(\Omega) = \intO\skp*{\nu_0;s}\dx{x} + m_0^\vr(\Omega),
        \]
        for any time $\tau>0$. Therefore $M_{\bar{\vr}}$ is a constant of motion also for the measure-valued solutions.

        \item The measure-valued solutions generated by an approximation as in Section~\ref{sec:existence} enjoy property~\eqref{eq:gradientmeasure}, since the Young measure is generated (in its third coordinate) by a sequence of gradients; and would still be true for other conceiveable approximations.         
This fact will be used to close the relative energy estimate. Note however that condition~\eqref{eq:gradientmeasure} itself does not imply that the projection $\pi_3\nu$ of the Young measure onto the third coordinate is generated by a sequence of gradients.
    \end{enumerate}
\end{remark}

\begin{remark}[Relation between two formulations]
The above definition of a measure-valued solution is formulated for the system~\eqref{eq:euler_2}. In particular we take advantage of formula~\eqref{eq:new_form} to obtain a conservative form of the nonlocal term. In order to relate this definition to the initial formulation~\eqref{eq:euler_1}, we need to make sense of the term $\rho\nabla_x\Phi_\vr$. This can be done by defining a distribution $\stackrel{\photon}{\rho\nabla_x\Phi_\rho}$ via the relation
\begin{equation*}
    \brk*{\stackrel{\photon}{\rho\nabla_x\Phi_\rho};\phi} = \int_0^\tau\int_\Omega\frac12\ol{|\nabla_x\Phi_\rho|^2}\divx\phi \dx{x}\dx{t} - \int_0^\tau\int_\Omega\ol{\nabla_x\Phi_\rho\otimes\nabla_x\Phi_\rho}:\nabla_x\phi\dx{x}\dx{t},
\end{equation*}
for any $\phi\in C^1(\bar\Omega)$, where $\brk*{\cdot;\cdot}$ denotes the duality pair between $(C^1)^{*}$ and $C^{1}$. It can be seen that the approximate sequences constructed in the next section give rise to such an object.
\end{remark}


\section{Existence of dissipative measure-valued solutions}
\label{sec:existence}
In this section we will prove existence of dissipative measure-valued solutions defined in the previous section.
To this end we construct a two-step approximation. Firstly, we define an approximate problem to~\eqref{eq:euler_1}, where we introduce a sixth order differential operator into the momentum equation. Secondly, we prove existence of a finite-dimensional approximation, and pass to the limit to obtain existence of weak solutions for the approximate problem. Finally, we pass to the limit with the coefficient of the highest order term to show that a sequence of such solutions generates a dissipative measure-valued solution to~\eqref{eq:euler_1}. We only consider the case when the initial data is regular, \ie, the initial concentration measures are all zero and $\nu_0$ is a Dirac measure concentrated at some given measurable initial states.
The main result of this section is the following.

\begin{theorem}\label{thm:existenceDMVS}
Let $1\leq d\leq 3$. Suppose $\Omega\subset\R^d$ is a bounded smooth domain. 
If the initial data $(\vr_0, \vu_0)$ is such that $\rho_0, \rho_0\vu_0\in L^1(\Omega)$ and has finite energy, then there exists a dissipative measure-valued solution with initial data 
\begin{equation*}
    \nu_{0,x} = \delta_{\{\rho_0(x), \vu_0(x), \nabla_x\Phi_{\rho}(0,x)\}}\;\;\;\text{ for a.e.\ }x\in\Omega.
\end{equation*}
\end{theorem}

\noindent The remainder of this section is dedicated to the proof of the above theorem.

\medskip

\noindent\emph{Approximate problem.} Let $\ep>0$ and let $(\!(\cdot;\cdot)\!) = (\cdot;\cdot)_{W_0^{3,2}(\Omega)}$ denote the standard scalar product in $W_0^{3,2}(\Omega)$.
Let $\rho_{0,\ep}$ and $\vu_{0,\ep}$ denote smooth functions obtained by standard mollification of $\rho_0$ and $\vu_0$ at scale $\ep$.
We say that a triple $(\rho^\ep, \vu^\ep, \Phi_{\rho^\ep})$ is a weak solution to the approximate Euler-Poisson problem with initial data $\rho_0^\ep = \rho_{0,\ep}+\ep$ and $\vu_0^\ep =\vu_{0,\ep}$ if, for all $\tau\in(0,T]$,

\begin{equation}\label{eq:approxcontinuity}
    \int_0^\tau\intO \rho^\ep \partial_t\psi + \rho^\ep \vu^\ep\cdot\nabla_x\psi\dx{x}\dx{t} = \intO\rho^\ep(\tau,x)\psi(\tau,x)\dx{x} - \intO\rho_0^\ep\psi(0,x)\dx{x},
\end{equation}
for all $\psi\in C^1([0,T]\times\overline{\Omega})$,
\begin{equation}\label{eq:approxmomentum}
\begin{split}
    \int_0^\tau\intO\rho^\ep \vu^\ep\cdot\partial_t\phi &+ \brk*{\rho^\ep \vu^\ep\otimes \vu^\ep - \nabla_x\Phi_{\rho^\ep}\otimes\nabla_x\Phi_{\rho^\ep}+\frac12|\nabla_x\Phi_{\rho^\ep}|^2\Id}:\nabla_x\phi\dx{x}\dx{t} -\gamma \int_0^\tau\intO \rho^\ep\vu^\ep\cdot\phi\dx{x}\dx{t} 
    \\[0.5em]
    &\qquad-\int_0^\tau\intO \rho^\ep x\cdot\phi\dx{x}\dx{t} - M_{\vr^\ep} \int_0^\tau\intO \nabla_x\Phi_{\vr^\ep}\cdot\phi \dx{x}\dx{t}
    \\[0.5em]
    &= \ep\int_0^\tau(\!(\vu^\ep;\phi)\!)\dx{t} + \intO(\rho^\ep\vu^\ep)(\tau,x)\cdot\phi(\tau,x)\dx{x} - \intO \rho_0^\ep \vu_0^\ep \cdot\phi(0,x)\dx{x},
\end{split}
\end{equation}
for all $\phi\in C^1([0,T]\times\overline{\Omega})$,
and
\begin{equation}\label{eq:approxpoisson}
    \intO\nabla_x\Phi_{\rho^\ep}(\tau,x)\cdot\nabla_x\vartheta(x)\dx{x} = \intO\rho^\ep(\tau,x)\vartheta(x)\dx{x}  -M_{\vr^\ep} \intO \vartheta(x) \dx{x},
\end{equation}
for all $\vartheta\in C^1(\overline{\Omega})$, where $M_{\vr^\ep} = \intO \rho^\ep(t,x) \dx{x}$ is constant in time for each $\ep>0$. 
Let us note that due to strong $L^1$ convergence of $\rho_0^\ep$ to $\rho_0$, we obtain convergence of $M_{\rho^\ep}$ towards $M_{\rho}$.

\noindent Let us remark that in the above formulation of the approximate problem we ignore the precise form of the sixth order elliptic operator as well as the corresponding boundary conditions. These conditions vanish when we pass to the limit with $\ep$, which is the sole purpose of this section.

We will now define Galerkin approximate solutions to the approximate Euler-Poisson problem. For convenience we drop the index $\ep$.

\noindent Let $\{\omega_i\}$ be an orthonormal basis of $W_0^{3,2}(\Omega)$ solving the eigenvalue problem
\[
(\!(\omega_i;\cdot)\!) = \lambda_i(\omega_i;\cdot)_{L^2(\Omega)}
\]
in $W_0^{3,2}(\Omega)$, \cf\cite[Section~6.4]{JosefBook}. Then $\{\omega_i\}$ is an orthonormal basis for $L^2(\Omega)$ and if the boundary of $\Omega$ is smooth enough, then all $\omega_i$ can be assumed smooth.
Now put $\vu^n(t,x) = \sum\limits_{i=1}^{n}c^n_i(t)\omega_i(x)$.

The triple $(\rho^n, \vu^n, \Phi_{\rho^n})$ is called a solution to the Galerkin approximation for the approximate Euler-Poisson problem~\eqref{eq:approxcontinuity}--\eqref{eq:approxpoisson} in $(0,T)\times\Omega$  with initial data $\rho^n_0 = \rho_0$ and $\vu_0^n = \sum\limits_{i=1}^{n}(\vu_0;\omega_i)_{L^2(\Omega)}\omega_i$, if
\begin{equation}\label{eq:galerkinmass}
    \partial_t\rho^n + \divx(\rho^n \vu^n) = 0,
\end{equation}
\begin{equation}\label{eq:galerkinmomentum}
    \begin{split}
        \intO&\brk*{\rho^n \partial_t \vu^n + \rho^n\nabla_x \vu^n\vu^n - \divx\brk*{ \nabla_x\Phi_{\rho^n}\otimes\nabla_x\Phi_{\rho^n}-\frac12|\nabla_x\Phi_{\rho^n}|^2\Id} + \gamma\rho^n \vu^n + \rho^n x}\cdot\omega_i\dx{x}
        \\[0.5em]
        &\hspace{2.5cm} + M_{\vr^n}\intO \nabla_x\Phi_{\vr^n}\cdot\omega_i \dx{x}  + \ep(\!(\vu^n;\omega_i)\!) = 0,
    \end{split}
\end{equation}
for $i=1,\dots,n$, and 
\begin{equation}\label{eq:galerkinpoisson}
    -\Delta\Phi_{\rho^n} = \rho^n - M_{\vr^n},
\end{equation}
where $M_{\vr^n} = \intO\rho^n(t,x)\dx{x}$.

\medskip

\noindent\emph{Existence of Galerkin approximations.} 
Observe that the initial data for the Galerkin problem satisfies in particular
$\rho_0^n\in C^1(\overline{\Omega})$, $\rho^n_0>0$ and $\vu^n_0\in W_0^{3,2}(\Omega)$.
For data with this regularity, the existence of a solution $(\rho^n, \vu^n)$ with
\[
    \rho^n\in C^1([0,T)\times\overline{\Omega}),\;\;\; \vu^n\in C^1([0,T); W_0^{3,2}(\Omega))
\]
to the Galerkin system
is obtained via a combination of the method of characteristics (to find the unique $\rho^n$) and the Schauder fixed point theorem. Since the proof follows in essentially the same steps as in~\cite{JosefBook} (see also~\cite{Gw2005}), we skip the details here. We just note that solving the continuity equation~\eqref{eq:galerkinmass} with the method of characteristics yields the following representation
\begin{equation}\label{eq:galerkindensity}
    \rho^n(t,x) = \rho_0(X^n(0;t,x))\exp{\prt*{-\int_0^t\divx \vu^n(\tau,X^n(\tau;t,x))\dx{\tau}}},
\end{equation}
where $X^n$ denotes the forward flow assosciated with velocity $\vu^n$.
In particular $\rho^n\geq \rho^*$ for some constant $\rho^*$ depending on $\ep^{-1}$.

\medskip

\noindent\emph{Energy estimates for the Galerkin approximations.}
Upon multiplying each equation of~\eqref{eq:galerkinmomentum} with the coefficient $c_i^n(t) = (\vu^n(t,\cdot);\omega_i(\cdot))_{L^2(\Omega)}$ and summing over $i=1,\dots,n$, we get, for each time,
\begin{equation*}
    \begin{aligned}
    \ddt\intO\frac12\rho^n|\vu^n|^2\dx{x} - \intO&\divx\brk*{\nabla_x\Phi_{\rho^n}\otimes\nabla_x\Phi_{\rho^n}-\frac12|\nabla_x\Phi_{\rho^n}|^2\Id} \cdot \vu^n\dx{x} + \ep(\!(\vu^n;\vu^n)\!)
    \\[0.5em]
    &= \gamma\intO \rho^n \abs*{\vu^n}^2 \dx{x} - \intO \rho^n x \cdot \vu^n\dx{x} - M_{\vr^n}\intO \nabla_x\Phi_{\vr^n}\cdot\vu^n \dx{x},
    \end{aligned}
\end{equation*}
where we used~\eqref{eq:galerkinmass} to write
\begin{equation*}
    \intO\rho^n \vu^n\cdot\nabla_x \vu^n \vu^n\dx{x} = \intO\partial_t\rho^n\ \frac12|\vu^n|^2\dx{x}.
\end{equation*}
Moreover, using~\eqref{eq:galerkinmass},~\eqref{eq:new_form} and the Poisson equation we have
\begin{equation*}
\begin{aligned}
    \intO\divx\brk*{\nabla_x\Phi_{\rho^n}\otimes\nabla_x\Phi_{\rho^n}-\frac12|\nabla_x\Phi_{\rho^n}|^2\Id}\cdot \vu^n\dx{x} - M_{\vr^n}\intO \nabla_x\Phi_{\vr^n}\cdot\vu^n \dx{x}&= -\intO\rho^n\nabla_x\Phi_{\rho^n}\cdot \vu^n\dx{x}
    \\[0.5em]
    &= -\ddt\intO\frac12|\nabla_x\Phi_{\rho^n}|^2\dx{x}.
\end{aligned}
\end{equation*}
Finally, we write
\begin{equation*}
    -\intO \rho^n x \cdot \vu^n \dx{x} = -\intO \rho^n\vu^n \cdot \nabla_x\prt*{\frac12\abs*{x}^2} \dx{x} =  -\ddt\intO\frac12\rho^n|x|^2\dx{x}.
\end{equation*}
We therefore obtain the following energy estimate
\begin{equation}\label{eq:galerkinenergy}
    \begin{split}
        \intO&\left[\frac12\rho^n |\vu^n|^2 + \frac12|\nabla_x\Phi_{\rho^n}|^2 + \frac12\rho^n|x|^2\right]\prt*{\tau,x}\dx{x} + \ep\int_0^\tau(\!(\vu^n;\vu^n)\!)\dx{t} 
        \\[0.5em]
        &\leq \intO\frac12\rho_0|\vu_0|^2 + \frac12|\nabla_x\Phi_{\rho^n}|^2(0,x) + \frac12\rho_0|x|^2\dx{x} - \gamma \int_0^\tau\intO\rho^n|\vu^n|^2\dx{x}\dx{t}. 
    \end{split}
\end{equation}
Now using~\eqref{eq:galerkindensity} and~\eqref{eq:galerkinenergy} we can deduce the uniform (in $n$) estimate
\begin{equation}\label{eq:galerkindensityestimate}
    \norm*{\rho^n}_{L^\infty((0,T)\times\Omega)} + \int_0^T\norm*{\partial_t\rho^n}^2_{L^2(\Omega)}\dx{t} + \int_0^T\norm*{\nabla_x\rho^n}^2_{L^2(\Omega)}\dx{t} \leq C\prt*{\ep^{-1}},
\end{equation}
while using~\eqref{eq:galerkinmomentum} multiplied by $\partial_t c_i^n(t)$,~\eqref{eq:galerkinenergy} and~\eqref{eq:galerkindensityestimate} we obtain
\begin{equation}\label{eq:galerkinvelocityestimate}
    \int_0^T\norm*{\partial_t \vu^n}^2_{L^2(\Omega)}\dx{t} + \ep\norm*{\vu^n}_{L^\infty\prt*{0,T;W_0^{3,2}(\Omega)}} \leq C\prt*{\ep^{-1}},
\end{equation}
see~\cite[Lemma~5.45]{JosefBook} for details.

\medskip

\noindent\emph{Existence of solutions to the approximate Euler-Poisson problem.}
Having~\eqref{eq:galerkinenergy}--\eqref{eq:galerkinvelocityestimate} and the Poisson equation~\eqref{eq:galerkinpoisson} we can, for each $\ep>0$, deduce, up to extracting a subsequence, the convergences

\begin{equation*}
    \begin{aligned}
        \rho^n &\weakstar \rho \;\;\; &&\text{ in } L^\infty((0,T)\times\Omega),\\
        \partial_t\rho^n &\weakto \partial_t\rho \;\;\; &&\text{ in } L^2((0,T)\times\Omega),\\
        \partial_t \vu^n &\weakto \partial_t \vu \;\;\; &&\text{ in } L^2((0,T)\times\Omega),\\
        \vu^n &\weakto \vu \;\;\; &&\text{ in } L^2((0,T);W_0^{3,2}(\Omega)),\\
        \nabla_x\Phi_{\rho^n}&\to\nabla_x\Phi_{\rho}\;\;\; &&\text{ in } L^2((0,T)\times\Omega).
    \end{aligned}
\end{equation*}
By the Aubin-Lions Lemma we then have
\begin{equation*}
    \begin{aligned}
    \rho^n &\to \rho \;\;\; &&\text{ in } L^2((0,T)\times\Omega),\\
    \vu^n &\to \vu \;\;\; &&\text{ in } L^2((0,T);W_0^{1,2}(\Omega)).
    \end{aligned}    
\end{equation*}
Combining the above weak and strong convergences, we can pass to the limit in each integral in the formulation of the Galerkin problem, thus showing existence of a solution $(\rho^\ep,\vu^\ep,\Phi_\rho^\ep)$ to the approximate Euler-Poisson problem.

\medskip

\noindent\emph{Existence of dissipative measure-valued solutions.} 
Since the approximate solution $(\rho^\ep, \vu^\ep, \Phi_\rho^\ep)$ is the limit of the Galerkin approximations, we have the energy bound
\begin{equation}\label{eq:approximateenergy}
    \begin{split}
        \intO&\left[\frac12\rho^\ep |\vu^\ep|^2 + \frac12|\nabla_x\Phi_{\rho^\ep}|^2 + \frac12\rho^\ep|x|^2\right]\prt{\tau,x}\dx{x} + \gamma \int_0^\tau\intO\rho^\ep|\vu^\ep|^2\dx{x}\dx{t}
        \\[0.5em]
        &\quad\leq \intO\left[\frac12\rho_0|\vu_0|^2 + \frac12|\nabla_x\Phi_{\rho^\ep}|^2(0,x) + \frac12\rho_0|x|^2\right]\dx{x} . 
    \end{split}
\end{equation}
Moreover, mass conservation implies that $\rho^\ep$ is uniformly bounded in $L^\infty(0,T;L^1(\Omega))$. Then,
\begin{equation*}
    \intO|\rho^\ep \vu^\ep|\dx{x}\leq\frac12\intO\rho^\ep\dx{x} + \frac12\intO\rho^\ep|\vu^\ep|^2\dx{x}\leq C.
\end{equation*}
Therefore the sequence of approximate momenta $\{\rho^\ep\vu^\ep\}$ is uniformly bounded in $L^\infty(0,T;L^1(\Omega))$, while $\{\nabla_x\Phi_{\rho^\ep}\}$ is bounded in $L^\infty(0,T;L^2(\Omega))$.

As discussed in the introduction, by considering the sequence $\delta_{\{\rho^\ep, \vu^\ep, \nabla_x\Phi_{\rho^\ep}\}}$ we obtain in the limit $\ep\to0$ a parameterised measure
\begin{equation*}
	\vnu=\{\nu_{t,x}\}\in L_{\weak}^\infty\prt*{(0,T)\times\Omega;\mathcal{M}^+\prt*{[0,\infty)\times		\R^d\times\R^d}}
\end{equation*}
which represents weak-$^*$ limits of nonlinear compositions with $C_0$ nonlinearities. However, since the ones appearing in the weak formulation of our problem are only continuous and not $C_0$  (and since their compositions with the approximating sequence are not uniformly integrable), we cannot apply any form of the Fundamental Theorem directly. The only terms in which we can straightaway pass to the limit are the ones containing $\nabla_x\Phi_{\rho^\ep}$ in the momentum and Poisson equations. For all the other terms we have to make use of Lemma~\ref{lem:integrability} to describe the oscillatory behaviour and introduce the concentration-defect measures as in~\eqref{eq:concentrationdefect}. Indeed the functions of interest are $f(s,\vv,\vf) = s, s\vv, s\vv\otimes\vv, s|\vv|^2, \vf\otimes\vf, |\vf|^2$, and all satisfy conditions of the lemma.

In particular, using~\eqref{eq:limitwithconcentration}, we can pass to the limit $\ep\to0$ in each term of~\eqref{eq:approxcontinuity}:
\begin{equation*}
    \begin{aligned}
    \int_0^\tau\intO \rho^\ep \partial_t\psi + \rho^\ep \vu^\ep\cdot\nabla_x\psi\dx{x}\dx{t} &\longrightarrow \int_0^\tau\intO \ol{\rho}\partial_t\psi + \ol{\rho\vu}\cdot\nabla_x\psi\dx{x}\dx{t},\\[0.5em]
    \intO\rho^\ep(\tau,x)\psi(\tau,x)\dx{x} &\longrightarrow \intO \ol{\rho}(\tau)\psi(\tau,x)\dx{x}, 
    \end{aligned}
\end{equation*}
for any $\psi\in C^1([0,T]\times\ol{\Omega})$, to obtain~\eqref{eq:mvcontinuity}.
Similarly we pass to the limit in~\eqref{eq:approxpoisson} to get~\eqref{eq:mvPoisson} (notice that $m^{\nabla\Phi}=0$, since $\nabla_x\Phi_\ep$ is square-integrable), and we also obtain~\eqref{eq:mvmomentum} as a limit of~\eqref{eq:approxmomentum}.
Finally, passing to the limit in~\eqref{eq:approximateenergy} we obtain~\eqref{eq:energy_inequality}. 
Notice that by virtue of Proposition~\ref{prop:concentrationrelations} we have several natural relations between the concentration measures, for instance
\begin{equation*}
	\abs*{m^{\rho|\vu|}} \leq m^{\rho} + m^{\rho|\vu|^2},
\end{equation*}
and
\begin{equation*}
	\abs*{m^{\nabla\Phi\otimes\nabla\Phi}} \leq m^{|\nabla\Phi|^2},\quad \abs*{m^{\rho\vu\otimes\vu}} \leq m^{\rho|\vu|^2}.
\end{equation*}
In particular $m^{|\nabla\Phi|^2} = 0$ implies $m^{\nabla\Phi\otimes\nabla\Phi} = 0$ (see Lemma~\ref{lem:inegalites2}).
Furthermore, due to bounds in $L^\infty(0,T;L^(\Omega))$, each concentration measure admits a disintegration with respect to the time and space variable, and thus we have
\begin{align*}
    m^\rho, m^{\rho|\vu|^2}, m^{|\nabla\Phi|^2} &\in L^\infty\prt*{0,T;\M^+(\ol{\Omega})},\\
    m^{\rho\vu} &\in L^\infty\prt*{0,T;\M(\ol{\Omega})^d},\\
    m^{\rho\vu\otimes\vu}, m^{\nabla\Phi\otimes\nabla\Phi} &\in L^\infty\prt*{0,T;\M(\ol{\Omega})^{d\times d}},
\end{align*}
Thus the proof of Theorem~\ref{thm:existenceDMVS} is complete.

\bigskip
Let us point out one useful identity satisfied by the measure-valued solutions. Suppose $(r,\vU, \Phi_r)$ is a regular solution of~\eqref{eq:euler_2}. Upon using~\eqref{eq:new_form} for this regular solution and for the Galerkin approximates, we obtain, in the limit as $n\to\infty$, the following identity

\begin{equation*}
\begin{aligned}
    &\int_0^\tau\intO \prt*{r-M_r}\nabla_x\Phi_{\rho^\ep}\cdot\phi \dx{x}\dx{t} + \int_0^\tau\intO \prt*{\rho^\ep-M_{\vr^\ep}}\nabla_x\Phi_r\cdot\phi \dx{x}\dx{t}
    \\
    &= -\int_0^\tau\intO \nabla_x\Phi_{\rho^\ep}\cdot\nabla_x\Phi_r\divx\phi \dx{x}\dx{t}
    \\
    &\quad + \int_0^\tau\intO \nabla_x\Phi_{\rho^\ep}\otimes\nabla_x\Phi_r\cdot\nabla_x\phi \dx{x}\dx{t} + \int_0^\tau\intO \nabla_x\Phi_r\otimes\nabla_x\Phi_{\rho^\ep}\cdot\nabla_x\phi \dx{x}\dx{t},
\end{aligned}
\end{equation*}
for every $\phi\in L^\infty(0,T;C^{1}(\Omega))$.

We can pass to the limit in each term of the above identity to obtain
\begin{equation}
\begin{aligned}\label{eq:mvnewform2}
    &\int_0^\tau\intO \prt*{r-M_r}\ol{\nabla_x\Phi_\rho}\cdot\phi \dx{x}\dx{t} + \int_0^\tau\intO \ol{\prt*{\rho-M_\vr}}\nabla_x\Phi_r\cdot\phi \dx{x}\dx{t}
    \\
    &= -\int_0^\tau\intO \ol{\nabla_x\Phi_\rho}\cdot\nabla_x\Phi_r\divx\phi \dx{x}\dx{t}
    \\
    &\quad + \int_0^\tau\intO \ol{\nabla_x\Phi_\rho\otimes\nabla_x\Phi_r}\cdot\nabla_x\phi \dx{x}\dx{t} + \int_0^\tau\intO \ol{\nabla_x\Phi_r\otimes\nabla_x\Phi_\rho}\cdot\nabla_x\phi \dx{x}\dx{t}.
\end{aligned}
\end{equation}

\begin{remark}[Other approximation schemes]
 Admittedly, the above approximation scheme carries the disadvantage of imposing rather strange and nonphysical boundary condition in the velocity. To avoid this, other approximations are possible, for instance through the Navier-Stokes equations with artificial local pressure. We skip the full exposition here, and refer to \cite{FGSW2016} for some details.
\end{remark}

\begin{remark}
\label{rem:othervariables}
	Notice that in the above proof we can only guarantee that the measures $\nu_{t,x}$ belong to the space $\M^+$ of non-negative Radon measures with $\norm{\nu_{t,x}}_{TV}\leq 1$. In the variables $(\rho,\vu,\nabla_x\Phi_\rho)$ one cannot guarantee that thease are necessariliy probability measures, because the tightness condition might fail. Indeed, when $\rho=0$, neither the momentum nor the kinetic energy offer control over the velocity, which can be arbitrary in the vaccum regions. This situation is reminiscent of other Euler-type models analysed in the past: the Savage-Hutter equations~\cite{Gw2005} or compressible Euler equations~\cite{GSW2015}. There, to circumvent this issue a different set of variables is considered, namely $\rho$ and $\sqrt{\rho}\vu$, instead of the more traditional $\rho$ and $\vu$. In this formulation one has the uniform bound
\begin{equation*}
	\int_{B_R(0)} |\sqrt{\rho^\ep}\vu^\ep|^2 \dx{x} = \int_{B_R(0)} \rho^\ep|\vu^\ep|^2 \dx{x} < \infty,
\end{equation*}
which implies tightness. Therefore the sequence $(\rho^\ep,\sqrt{\rho^n}\vu^n,\nabla_x\Phi_\rho^n)$ generates a Young measure 
\[
	\vec{\mu} \in L^\infty_{\weak}\prt*{(0,T)\times\R^d ; \PP([0,\infty)\times\R^d\times\R^d)},
\]
which clearly agrees with $\vnu$ on the first and third coordinates.
\end{remark}


\section{Relative energy inequality}
\label{sec:relative}

We suppose now that $(r,\vU,\Phi_r)$, $r>0$, is a strong solution to the Euler-Poisson system~\eqref{eq:euler_1} with regular initial data $(\rho_0, \vu_0)$ of finite energy. 
Furthermore we consider a dissipative measure-valued solution $(\vnu, m^\rho, m^{\rho\vu}, m^{\rho \vu\otimes \vu}, m^{|\nabla\Phi|^2}, m^{\nabla\Phi\otimes\nabla\Phi})$ with
\[
    \nu_{0,x} = \delta_{\{\rho_0(x), \vu_0(x), \nabla_x\Phi_r(0,x)\}}
\]
for a.e.\ $x\in\Omega$.

In~\cite{CaFeGwSw2015} the following relative energy was used:
\[
    \E (\vr,\vu \,|\, r, \vU) = \int_{\Omega} 
    \brk*{
    \frac{1}{2} \vr \abs*{\vu - \vU}^2 + \frac{1}{2} (r - \vr)(K \ast (r-\vr))
    } \dx{x},
\]
$K$ being the Poisson kernel, to compare strong and dissipative weak solutions and establish a weak-strong uniqueness result. Note that in their case $\Omega = \mathbb{T}^d$ was the flat torus in two or three dimensions.
We shall mimic this approach here. To this end, we notice that, upon an integration by parts, we have
\[
    \intO \frac{1}{2} (r - \vr)(K \ast (r-\vr)) \dx{x} = \intO\frac12 \abs*{\nabla_x\Phi_r-\nabla_x\Phi_\rho}^2 \dx{x}.
\]
A natural candidate for a measure-valued version of the relative energy is therefore:
\begin{equation}\label{eq:mvrelenergy}
    \E_{rel}^{mv}(\tau)= \int_{\Omega} 
    \left[
    \frac{1}{2} \ol{\rho |\vu - \vU|^2} + \frac{1}{2} \ol{|\nabla_x\Phi_\rho - \nabla_x\Phi_r|^2}
    \right] \dx{x}.
\end{equation}

\noindent We can then write
\begin{equation}\label{eq:mvrelenergy2}
    \begin{split}
    \E_{rel}^{mv}(\tau) =& \intO\frac12 \ol{\rho |\vu|^2} \dx{x} + \intO\frac12 \ol{|\nabla_x\Phi_\rho|^2} \dx{x} + \intO\frac12 \ol{\rho}|\vU|^2\dx{x} 
    \\[0.5em]
    &- \intO\ol{\rho\vu}\cdot\vU\dx{x} - \intO \ol{\nabla_x\Phi_\rho}\cdot\nabla_x\Phi_r \dx{x} + \frac12\intO|\nabla_x\Phi_r|^2\dx{x}. 
    \end{split}
\end{equation}
We also introduce the measure-valued variant of the energy
\[
    \E^{mv}(\tau)= \int_{\Omega} 
    \left[
    \frac{1}{2} \ol{\rho |\vu|^2} + \frac{1}{2} \ol{|\nabla_x\Phi_\rho|^2} + \frac12 \ol{\rho}|x|^2
    \right] \dx{x},
\]
so that inequality~\eqref{eq:energy_inequality} becomes
\begin{equation}\label{eq:mvenergyinequality}
    \E^{mv}(\tau) \leq \E^{mv}(0) - \gamma \int_0^\tau \int_{\Omega} 
	\ol{\rho |\vu|^2} \dx{x}\dx{t} .
\end{equation}

\noindent Testing the continuity equation~\eqref{eq:mvcontinuity} in turn with $\frac12|\vU|^2$ and $\Phi_r$ we have
\begin{equation*}
    \begin{split}
    \intO\frac12\ol{\rho}|\vU|^2(\tau,x)\dx{x} = \intO&\frac12\rho_0|\vU_0|^2 + \int_0^\tau\intO\ol{\rho}\vU\cdot\partial_t\vU + \ol{\rho\vu}\cdot\nabla_x\vU\vU\dx{x}\dx{t},
    \end{split}
\end{equation*}
and
\begin{equation*}
    \begin{split}
    \intO\ol{\rho}\Phi_r(\tau,x)\dx{x} = \intO&\rho_0\Phi_r(0,x) + \int_0^\tau\intO\ol{\rho}\partial_t\Phi_r + \ol{\rho\vu}\cdot\nabla_x\Phi_r\dx{x}\dx{t},
    \end{split}
\end{equation*}
while testing the momentum equation~\eqref{eq:mvmomentum} with $\vU$ gives
\begin{equation*}
\begin{split}
	\int_{\Omega} &\ol{\rho \vu}\cdot \vU (\tau, x) \dx{x}
	= \int_{\Omega} \rho_0|\vU_0|^2 \dx{x} 
	+  \int_0^\tau \int_\Omega 
	\left[
	\ol{\rho \vu}\cdot \partial_t \vU 
	+ \ol{\rho \vu\otimes\vu} : \nabla_x \vU
	- \gamma \ol{\rho \vu}\cdot \vU
	- \ol{\rho} x\cdot\vU
	\right] \dx{x}\dx{t}
	\\[0.5em]
	&\hspace{3cm}+\int_0^\tau \int_\Omega \frac12\ol{|\nabla_x\Phi_\rho|^2}\, \divx\vU -\ol{\nabla_x\Phi_\rho\otimes\nabla_x\Phi_\rho} : \nabla_x\vU - M_\vr \ol{\nabla_x\Phi_\vr}\cdot\vU \dx{x}\dx{t}.
\end{split}
\end{equation*}
Note also that by~\eqref{eq:mvPoisson} tested with $\Phi_r$ we have
\begin{equation*}
    \intO\ol{\rho}\Phi_r(\tau,x)\dx{x} - M_\vr\int\Phi_r(\tau,x)\dx{x} = \intO\ol{\nabla_x\Phi_\rho}\cdot\nabla_x\Phi_r(\tau,x)\dx{x}.
\end{equation*}

\noindent Using the above identities and~\eqref{eq:mvrelenergy2} we get

\begin{equation*}
    \begin{split}
    \E_{rel}^{mv}(\tau)
    &= \E^{mv}(\tau) - \intO\frac12\rho_0|\vU_0|^2\dx{x} - \intO\rho_0\Phi_r(0,x)\dx{x}
    \\[0.5em]
    &\quad+\int_0^\tau\intO\ol{\rho}\vU\cdot\partial_t\vU + \ol{\rho\vu}\cdot\nabla_x\vU\vU-\ol{\rho \vu}\cdot \partial_t \vU - \ol{\rho \vu\otimes\vu} : \nabla_x \vU\dx{x}\dx{t}
    \\[0.5em]
	&\quad+\gamma\int_0^\tau\intO \ol{\rho \vu}\cdot \vU \dx{x}\dx{t}
	+ \int_0^\tau\intO\ol{\rho} x\cdot\vU \dx{x}\dx{t}
	+\frac12\intO|\nabla\Phi_r|^2(\tau,x)\dx{x}
	\\[0.5em]
	&\quad 
	-\int_0^\tau\intO\frac12\ol{|\nabla_x\Phi_\rho|^2} \divx\vU - \ol{\nabla_x\Phi_\rho\otimes\nabla_x\Phi_\rho} : \nabla_x\vU - M_\vr \ol{\nabla_x\Phi_\vr}\cdot\vU \dx{x}\dx{t}
	\\[0.5em]
	&\quad-\int_0^\tau\intO\ol{\rho}\partial_t\Phi_r\dx{x}\dx{t} -\int_0^\tau\intO\ol{\rho\vu}\nabla_x\Phi_r\dx{x}\dx{t} - \frac12\intO\ol{\rho}(\tau,x)|x|^2\dx{x}
	\\[0.5em]
	&\quad + M_\vr\intO\Phi_r(\tau,x)\dx{x}.
    \end{split}
\end{equation*}
For the first line on the right-hand side we have, invoking~\eqref{eq:mvenergyinequality} and~\eqref{eq:mvPoisson} again, an upper bound given by
\[
    \frac12\intO\rho_0|x|^2\dx{x} - \frac12\intO|\nabla_x\Phi_r|^2(0,x)\dx{x}
    - \gamma\int_0^\tau \int_{\Omega}\ol{\rho |\vu|^2} \dx{t} -M_\vr\intO\Phi_r(0,x)\dx{x},
\]
while for the second line we write, using the momentum equation for the strong solution $(r,\vU)$ and the strict positivity of $r$,
\begin{equation*}
    \begin{split}
    \int_0^\tau\intO&\ol{\rho}\vU\partial_t\vU + \ol{\rho\vu}\cdot\nabla_x\vU\vU-\ol{\rho\vu}\cdot\partial_t \vU 
	- \ol{\rho\vu\otimes\vu} : \nabla_x \vU\dx{x}\dx{t}
	\\[0.5em]
	&= \int_0^\tau\intO\ol{\rho(\vU-\vu)}\cdot\partial_t\vU + \ol{\rho\vu\otimes(\vU-\vu)}:\nabla_x\vU\dx{x}\dx{t}
	\\[0.5em]
	&= \int_0^\tau\intO\ol{\rho(\vu-\vU)}\cdot(\gamma\vU + x + \nabla_x\Phi_r) + \ol{\rho(\vu-\vU)\otimes(\vU-\vu)}:\nabla_x\vU\dx{x}\dx{t}.
    \end{split}    
\end{equation*}
Furthermore we write
\begin{equation*}
    \begin{split}
    \int_0^\tau\intO\ol{\rho}x\cdot\vU\dx{x}\dx{t} &= \int_0^\tau\intO x\cdot r\vU\dx{x}\dx{t} + \int_0^\tau\intO\ol{(\rho-r)} x\cdot\vU\dx{x}\dx{t}
    \\[0.5em]
    &= \intO\frac12|x|^2(r(\tau,x)-r_0)\dx{x}+ \int_0^\tau\intO\ol{(\rho-r)} x\cdot\vU\dx{x}\dx{t}.
    \end{split}
\end{equation*}
We thus have
\begin{equation*}
    \begin{split}
    \E_{rel}^{mv}(\tau) \leq &\left[ \frac12\intO\ol{(r-\rho)}|x|^2 \dx{x}\right]_{t=0}^{t=\tau}
    +\int_0^\tau\intO\ol{\rho(\vu-\vU)\otimes(\vU-\vu)}:\nabla_x\vU\dx{x}\dx{t}
    \\[0.5em]
	&+\gamma\int_0^\tau\intO \ol{\rho \vu} \cdot \vU - \ol{\rho|\vu|^2} + \ol{\rho(\vu-\vU)}\cdot\vU\dx{x}\dx{t}
	\\[0.5em]
	&+\int_0^\tau\intO\ol{(\rho-r)} x\cdot\vU + \ol{\rho(\vu-\vU)}\cdot x\dx{x}\dx{t}
	\\[0.5em]
	&+\frac12\intO|\nabla_x\Phi_r|^2(\tau,x)\dx{x} - \frac12\intO|\nabla_x\Phi_r|^2(0,x)\dx{x} -\int_0^\tau\intO\ol{\rho}\partial_t\Phi_r\dx{x}\dx{t} +\left[ M_\vr\intO \Phi_r(t,x)\dx{x}\right]_{t=0}^{t=\tau}
	\\[0.5em]
	&-\int_0^\tau\intO \brk*{\frac12\ol{|\nabla_x\Phi_\rho|^2}\divx\vU -\ol{\nabla_x\Phi_\rho\otimes\nabla_x\Phi_\rho} : \nabla_x\vU - M_\vr \ol{\nabla_x\Phi_\vr}\cdot\vU + \ol{\rho}\vU\cdot\nabla_x\Phi_r} \dx{x}\dx{t}.
\end{split}
\end{equation*}

\noindent Notice that the function $(t,x)\mapsto\frac12|x|^2$ is an admissible test function for the continuity equation~\eqref{eq:mvcontinuity}, so that
\begin{equation*}
    \left[ \frac12\intO\ol{(r-\rho)}|x|^2\dx{x}\right]_{t=0}^{t=\tau} = -\int_0^\tau\intO\ol{(\rho\vu-r\vU)}x\dx{x}\dx{t}.
\end{equation*}
Furthermore, we make note of the following identities:
\begin{equation*}
    \begin{aligned}
    \gamma\int_0^\tau\intO \ol{\rho \vu} \cdot \vU - \ol{\rho|\vu|^2} + \ol{\rho(\vu-\vU)}\cdot\vU\dx{x}\dx{t} &= -\gamma\int_0^\tau\intO \ol{\rho |\vu-\vU|^2}\dx{x}\dx{t},
    \\[0.5em]
    \int_0^\tau\intO\ol{(\rho-r)} x\cdot\vU + \ol{\rho(\vu-\vU)}\cdot x\dx{x}\dx{t} &= \int_0^\tau\intO\ol{(\rho\vu-r\vU)}x\dx{x}\dx{t},
    \end{aligned}
\end{equation*}
and
\begin{equation*}
    \begin{aligned}
    \frac12\intO|\nabla_x\Phi_r|^2(\tau,x)\dx{x} &- \frac12\intO|\nabla_x\Phi_r|^2(0,x)\dx{x} -\int_0^\tau\intO\ol{\rho}\partial_t\Phi_r\dx{x}\dx{t} +\left[ M_\vr\intO \Phi_r(t,x)\dx{x}\right]_{t=0}^{t=\tau}
    \\[0.5em]
    &= \frac12\int_0^\tau\intO\partial_t|\nabla_x\Phi_r|^2\dx{x}\dx{t}-\int_0^\tau\intO\ol{\rho}\partial_t\Phi_r\dx{x}\dx{t}+M_\vr\int_0^\tau\intO \partial_t\Phi_r\dx{x}\dx{t}
    \\[0.5em]
    &=\int_0^\tau\intO\nabla_x\Phi_r\cdot\nabla_x(\partial_t\Phi_r)\dx{x}\dx{t}-\int_0^\tau\intO\ol{\rho}\partial_t\Phi_r\dx{x}\dx{t}+M_\vr\int_0^\tau\intO \partial_t\Phi_r\dx{x}\dx{t}
    \\[0.5em]
    &=\int_0^\tau\intO\ol{(\nabla_x\Phi_r-\nabla_x\Phi_\rho)}\nabla_x(\partial_t\Phi_r)\dx{x}\dx{t},
    \end{aligned}
\end{equation*}
where for the last equality we used~\eqref{eq:mvPoisson}.

\noindent Finally, we use~\eqref{eq:mvnewform2} to write
\begin{align*}
    &-\int_0^\tau\intO \brk*{\frac12\ol{|\nabla_x\Phi_\rho|^2}\divx\vU -\ol{\nabla_x\Phi_\rho\otimes\nabla_x\Phi_\rho} : \nabla_x\vU - M_\vr \ol{\nabla_x\Phi_\vr}\cdot\vU + \ol{\rho}\vU\cdot\nabla_x\Phi_r} \dx{x}\dx{t}
    \\[0.5em]
    &=-\frac12\int_0^\tau\intO\ol{|\nabla_x\Phi_\rho-\nabla_x\Phi_r|^2}\divx\vU \dx{x}\dx{t} + \int_0^\tau\intO\ol{(\nabla_x\Phi_\rho-\nabla_x\Phi_r)\otimes(\nabla_x\Phi_\rho-\nabla_x\Phi_r)}:\nabla_x\vU \dx{x}\dx{t}
    \\[0.5em]
    &\quad
    + \int_0^\tau\intO\ol{(\nabla_x\Phi_\rho-\nabla_x\Phi_r)}\cdot r\vU\dx{x}\dx{t}.
\end{align*}

\noindent Furthermore, we use the strong form of the Poisson equation to write
\begin{align*}
    \int_0^\tau&\intO\ol{(\nabla_x\Phi_\rho-\nabla_x\Phi_r)}\cdot r\vU\dx{x}\dx{t} +\int_0^\tau\intO\ol{(\nabla_x\Phi_r-\nabla_x\Phi_\rho)}\cdot\nabla_x(\partial_t\Phi_r)\dx{x}\dx{t}
    \\[0.5em]
    &=\int_0^\tau\intO\ol{\nabla_x\Phi_\rho}\cdot r\vU\dx{x}\dx{t} -\int_0^\tau\intO\ol{\nabla_x\Phi_\rho}\cdot\nabla_x(\partial_t\Phi_r)\dx{x}\dx{t}
    \\[0.5em]
    &=\int_0^\tau\intO\skp*{\nu_{t,x};\vf}\cdot r\vU\dx{x}\dx{t} -\int_0^\tau\intO\skp*{\nu_{t,x};\vf}\cdot\nabla_x(\partial_t\Phi_r)\dx{x}\dx{t}
    \\[0.5em]
    &=\int_0^\tau\intO\nabla_x\Psi_\rho(t,x)\cdot r\vU\dx{x}\dx{t} -\int_0^\tau\intO\nabla_x\Psi_\rho(t,x)\cdot\nabla_x(\partial_t\Phi_r)\dx{x}\dx{t}
    \\[0.5em]
    &=0,
\end{align*}
using property~\eqref{eq:gradientmeasure} in the third equality above.

\noindent We therefore arrive at the following measure-valued version of the relative energy inequality
\begin{equation}\label{eq:mvrelativeenergy}
    \begin{split}
    \E_{rel}^{mv}(\tau)\leq & \int_0^\tau\intO\ol{\rho(\vu-\vU)\otimes(\vU-\vu)}:\nabla_x\vU\dx{x}\dx{t}
    -\gamma\int_0^\tau\intO \ol{\rho |\vu-\vU|^2}\dx{x}\dx{t}
    \\[0.5em]
    &-\frac12\int_0^\tau\intO\ol{|\nabla_x\Phi_\rho-\nabla_x\Phi_r|^2}\divx\vU \dx{x}\dx{t} + \int_0^\tau\intO\ol{(\nabla_x\Phi_\rho-\nabla_x\Phi_r)\otimes(\nabla_x\Phi_\rho-\nabla_x\Phi_r)}:\nabla_x\vU \dx{x}\dx{t}.
    \end{split}
\end{equation}

\noindent All the terms on the right-hand side of inequality~\eqref{eq:mvrelativeenergy} can be easily seen to be controlled by 
\[
    c\int_0^\tau\E_{rel}^{mv}(t)\dx{t},
\]    
where the constants $c$ depend only on the norm $\norm{\nabla_x\vU}_{C([0,T]\times\ol{\Omega})}$. Note that for the terms involving a tensor product we use Lemma~\ref{lem:inegalites2} to compare them to the corresponding norm-squared terms.

We thus have the following Gronwall-type inequality
\begin{equation}\label{eq:relenergyfinal}
        \E_{rel}^{mv}(\tau) \leq c \int_0^\tau\E_{rel}^{mv}(t)\dx{t}.
\end{equation}

\section{Proof of the main theorem}
\label{sec:mainproof}

 Having established the relative energy inequality, we proceed to the proof of our main result, Theorem~\ref{thm:mv-stronguniqness}. This will be done in several steps, each leading to identifications of successive terms in the measure-valued formulation with their counterparts in the strong formulation.

Firstly, we observe that inequality~\eqref{eq:relenergyfinal} implies that $\E_{rel}^{mv} = 0$ at almost all times, since the strong and measure-valued solution emanate from the same initial data. In particular, since both terms of the relatie energy are non-negative, we have
\begin{equation}
    \intO\ol{|\nabla_x\Phi_\rho-\nabla_x\Phi_r|^2}\dx{x} = 0,
\end{equation}
from which we readily infer that the projection of the Young measure $\nu$ onto the third coordinate reduces to a Dirac mass at $\nabla_x\Phi_r$, and therefore
\begin{equation}
\label{eq:potentialidentification}
    \nu_{t,x}(s,\vv,\vf) = \bar{\nu}_{t,x}(s,\vv)\otimes\delta_{\{\nabla\Phi_r(t,x)\}},\quad m^{|\nabla\Phi|^2} = 0,
\end{equation}
where
\[
\bar{\nu} = \prt*{\bar{\nu}_{t,x}} \in L^\infty_{weak}\prt*{(0,T)\times\Omega ; \mathcal{M}^{+}([0,\infty)\times\R^d)}.
\]
Consequently, in the ``kinetic'' part of the relative energy we have
\begin{equation}
    \int_{[0,\infty)\times\R^d} s|\vv-\vU|^2\dx{\bar{\nu}_{t,x}}(s,\vv) = 0,\quad m^{\rho|\vu-\vU|^2} = 0.
\end{equation}
As explained in the introduction, however, we cannot conclude here the vanishing of the concentration measures in the density and the momentum. Instead, we need to work with whole oscillation-concentration pairs and only relate these to the corresponding strong quantities.

We begin by considering the kinetic term appearing in the energy inequality: namely, we will show that 
\[
\ol{\rho|\vu|^2} = \ol{\rho}|\vU|^2 = r|\vU|^2.
\]
To this end we choose $0<\delta<1$ and apply Lemma~\ref{lem:inegalites} from Appendix~\ref{sec:appB} to infer that on the level of approximating sequences $\rho^\ep, \vu^\ep$, as in Section~\ref{sec:existence}, we have the inequality
\begin{equation}
	\rho^\ep\abs*{|\vu^\ep|^2 - |\vU|^2} \leq C\delta\rho^\ep + C_\delta\rho^\ep|\vu^\ep-\vU|^2.
\end{equation}
Since these inequalities are preserved when passing with $n$ to infinity, we conclude
\begin{equation}
	-C\delta\ol{\rho} \leq \ol{\rho|\vu|^2} - \ol{\rho}|\vU|^2 \leq C\delta\ol{\rho},
\end{equation}
where we have used that
\begin{equation}
	\ol{\rho|\vu-\vU|^2} = 0,
\end{equation}
as concluded from the relative energy inequality.
Whence, by arbitrariness of $\delta>0$, we have
\begin{equation}
\label{eq:kineticidentification}
	\ol{\rho|\vu|^2} = \ol{\rho}|\vU|^2.
\end{equation}
Let us mention in passing that this equality concerns the sums of the Young measure and concentration parts of the corresponding terms and it is in principle not immediately obvious that corresponding term-by-term equalities follow (in particular that $m^{\rho|\vu|^2} = |\vU|^2m^\rho$). However, this is indeed true by virtue of Lemma~\ref{lem:projections}.\\
Next, we observe that
\begin{equation}
\label{eq:densityidentification}
	\ol{\rho} = r
\end{equation}
almost everywhere in $(0,T)\times\Omega$.
Indeed, this follows from the identification $\skp{\nu_{t,x};\vf} = \nabla_x\Phi_r(t,x)$, mass conservation and the Poisson equation. Therefore, we have
\begin{equation}
	\ol{\rho|\vu|^2} = r|\vU|^2,
\end{equation}
almost everywhere in $(0,T)\times\Omega$.

Similarly, applying Lemma~\ref{lem:inegalites} again, we have
\begin{equation}
	\rho^\ep\abs*{u^\ep_iu^\ep_j - U_iU_j} \leq C\delta\rho^\ep + C_\delta\rho^\ep|\vu^\ep-\vU|^2,
\end{equation}
and therefore in the limit
\begin{equation}
	\ol{\rho u_iu_j} = rU_iU_j,
\end{equation}
so that we can conclude the analogous equality for the convective term, \ie,
\begin{equation}
\label{eq:convectiveidentification}
	\ol{\rho\vu\otimes\vu} = r\vU\otimes\vU.
\end{equation}

Consequently, the momentum equation becomes the simple ODE
\begin{equation}
    \partial_t\prt*{\ol{\rho \vu}} + \gamma\ol{\rho\vu} = \partial_t\prt*{r\vU} + \gamma r\vu,
\end{equation}
from which we infer that
\begin{equation}
\label{eq:momentumidentification}
    \ol{\rho \vu}(t,x) = (r\vU)(t,x),
\end{equation}
for almost every $(t,x)\in(0,T)\times\Omega$.

Let us now focus on the Young measure $\bar{\nu}$ from equality~\eqref{eq:potentialidentification}. Firstly, to deal with potential "vacuum regions" where the density vanishes, we decompose this measure into 
\begin{equation}
	\bar{\nu} = \bar{\nu} \mres \prt*{\{0\}\times\R^d} + \bar{\nu} \mres \prt*{(0,\infty)\times\R^d} =\colon \sigma^1 + \sigma^2.
\end{equation}
Then, since
\begin{equation}
	\int_{[0,\infty)\times\R^d} s|\vv-\vU|^2 \dx{\bar{\nu}_{t,x}} = \int_{(0,\infty)\times\R^d} s|\vv-\vU|^2 \dx{\sigma^2_{t,x}} = 0, 
\end{equation}
we have
\begin{equation}
	\sigma^2_{t,x} = \bar{\bar{\nu}}_{t,x} \otimes \delta_{\{\vU(t,x)\}}
\end{equation}
almost everywhere in $(0,T)\times\Omega$ for some measure $\bar{\bar{\nu}}_{t,x}\in \mathcal{M}^+(0,\infty)$. Note that we still cannot guarantee that this is a probability measure.
Consequently, denoting by $\pi_1\sigma^1$ the projection of the measure $\sigma^1$ onto the first coordinate, and using Lemma~\ref{lem:projections} we may expand equality~\eqref{eq:densityidentification} as
\begin{equation}
	\int_A \brk*{\int_{\{0\}}s\dx{(\pi_1\sigma^1_{t,x})} + \int_{(0,\infty)}s\dx{\bar{\bar{\nu}}_{t,x}}} \dx{x}\dx{t} + m^\rho(A) = \int_A r(t,x) \dx{x}\dx{t},
\end{equation}
where $A$ is any Borel subset of $(0,T)\times\R^d$.
In particular, since $r>0$, this equality implies that the concentration measure $m^\rho$ is absolutely continuous with respect to the Lebesgue measure on $[0,\infty)\times\R^d$, so that there exists its Radon-Nikodym derivative, $D^{m^\rho}$. Hence, the last equality can be rewritten as
\begin{equation}
\label{eq:densityidentification2}
	\int_{(0,\infty)} s \dx{\bar{\bar{\nu}}_{t,x}}(s) + D^{m^\rho}(t,x) = r(t,x),
\end{equation}
almost everywhere.
In a similar manner equalities~\eqref{eq:kineticidentification},~\eqref{eq:momentumidentification},~\eqref{eq:convectiveidentification} imply that all the measure $m^{\rho|\vu|^2}$, $m^{\rho\vu}$, $m^{\rho\vu\otimes\vu}$ are absolutely continuous with respect to the Lebesgue measure.
This concludes the proof of Theorem~\ref{thm:mv-stronguniqness}.

\begin{remark}
	At this point one might make one more observation about the measure family of probabilities $\mu_{t,x}$ corresponding to the variables $(\rho, \sqrt{\rho}\vu, \nabla_x\Phi_\rho)$ as in Remark~\ref{rem:othervariables}. In this case one can deduce that
	\begin{equation}
		\mathrm{supp}\prt*{\pi_2\mu}_{t,x} \subset \{\vv\; :\; \vv = \alpha\vU(t,x), \alpha\in[0,\infty)\}
	\end{equation}	
so that the projection of the Young measure onto the second coordinate lives only on the one-dimensional subspace determined by the direction of the velocity $\vU$.		 
\end{remark}

\begin{remark}
Notice that if we consider the quantity
\begin{equation}
    \KE(t) \coloneqq \intO\frac12\ol{\rho|\vu|^2}\dx{x}.
\end{equation}
We then have
\begin{equation}
    \KE(\tau) \leq \E^{mv}(\tau) \leq \E^{mv}(0) - 2\int_0^\tau \KE(t)\dx{t}.
\end{equation}
Consequently $\KE$ satisfies, at almost all times, the bound
\begin{equation}
    0\leq \KE(\tau) \leq \E^{mv}(0)e^{-2t},
\end{equation}
and therefore converges to zero with time. Thus both the Young and the concentration measures in $\KE$ converge to zero as $t\to\infty$.
\end{remark}


\section{Euler alignment system}

The calculations performed in the previous section do not saturate the capacity of the relative energy method in the sense that we can also consider, instead of system~\eqref{eq:euler_1}, a system with more general form of the confinement potential, as well as include nonlinear damping.
More precisely, let us consider the following system Euler alignment system

\begin{equation}
	\begin{aligned}
	\partial_t \vr  + \divx (\vr \vu) & = 0,  
	\\
	\partial_t (\vr \vu) + \divx (\vr \vu \otimes \vu) & = - \vr \nabla_x \Phi_\vr - \vr \nabla_x V - \vr\nabla_x\prt*{W\star\rho} + \rho \intO\psi(x-y)\rho(t,y)\prt*{\vu(t,y)-\vu(t,x)}\dx{y},  
	\label{eq:euler_general}
	\\
	-\Delta \Phi_\vr &= \vr-M_\vr,  
	\end{aligned}
\end{equation}
in $(0,T) \times \Omega$, where $V=V(x)$ is smooth, and $W$ and $\psi$ are smooth and symmetric. The kernel $W$ includes the repulsive-attractive interaction force between individuals, while $\psi$ gives the local averaging measuring the consensus in their orientation. The equations are supplemented with the same boundary conditions~\eqref{eq:bdryconditions} as before. The spatial domain, $\Omega$, is still a bounded smooth domain in $\R^d$. Nevertheless, by a slight abuse of notation, we shall use the convolution symbol $\star$ to denote the integral
\begin{equation}
    (W\star\rho)(x) = \intO W(x-y)\rho(y)\dx{y},
\end{equation}
so as to simplify the notation.
A similar system has recently been considered in~\cite{BrezinaMacha}, where measure-valued solutions of a viscous approximation are shown to converge in the inviscid limit to a strong solution of equation~\eqref{eq:euler_general}. 
However, the presence of an artificial pressure term is assumed -- this substantially simplifies the analysis, as explained above.

The energy identity for strong solutions of system~\eqref{eq:euler_general} reads
\begin{equation}
    \begin{aligned}
    \ddt \intO \frac12\rho|\vu|^2 + \frac12\abs*{\nabla_x\Phi_\vr}^2 + \rho V + \frac12 \rho\prt*{W\star\rho} \dx{x}
    =
    -\frac12\intO\intO\psi(x-y)\rho(t,x)\rho(t,y)\abs*{\vu(t,y)-\vu(t,x)}^2 \dx{x}\dx{y}.
    \end{aligned}
\end{equation}

\noindent The definition of a dissipative measure-valued solutions has to be adjusted accordingly to the appearance of new terms. The weak formulation of the momentum equation, Equation~\eqref{eq:mvmomentum}, becomes

	\begin{equation}\label{eq:mvmomentum2}
	\begin{split}
	    \int_{\Omega} & \ol{\rho\vu}(\tau,x) \cdot \phi (\tau, x) \dx{x} 
	    - \int_{\Omega} \ol{\rho\vu}(0,x) \cdot \phi (0, x) \dx{x} 
	    \\[0.5em] 
	    & =
	    \int_0^\tau\!\! \int_\Omega \ol{\rho\vu}\cdot\partial_t \phi \dx{x}\dx{t}
	    +\int_0^\tau\!\! \int_\Omega \ol{\rho\vu\otimes\vu} : \nabla_x \phi \dx{x}\dx{t}
	    - \int_0^\tau\!\! \int_\Omega\ol{\rho}\, \nabla_x V\cdot\phi\dx{x}\dx{t}
	    \\[0.5em]
	    &\hspace{0.2cm}
	    + \int_0^\tau\!\! \int_\Omega\frac12\ol{|\nabla_x\Phi_\rho|^2}\,\divx\phi \dx{x}\dx{t} 
	    - \int_0^\tau\!\!\int_\Omega \ol{\nabla_x\Phi_\rho\otimes\nabla_x\Phi_\rho}:\nabla_x\phi \dx{x}\dx{t} 
	    - M_\rho\int_0^\tau\!\!\intO\ol{\nabla_x\Phi_\rho}\cdot\phi \dx{x} \dx{t}
	    \\[0.5em]
	    &\hspace{0.2cm}
	    - \int_0^\tau\!\!\intO \ol{\rho}(t,x)\intO\nabla_xW(x-y)\ol{\rho}(t,y)\dx{y}\cdot\phi(t,x)\dx{x}\dx{t}
	    \\[0.5em]
	    &\hspace{0.2cm}
	    + \int_0^\tau\!\!\intO \ol{\rho}(t,x)\intO\psi(x-y)\ol{\rho\vu}(t,y)\dx{y}\cdot\phi(t,x)\dx{x}\dx{t}
	    \\[0.5em]
	    &\hspace{0.2cm}
	    - \int_0^\tau\!\!\intO \ol{\rho\vu}(t,x)\intO\psi(x-y)\ol{\rho}(t,y)\dx{y}\cdot\phi(t,x)\dx{x}\dx{t}
	\end{split}
	\end{equation}
for a.a.\ $\tau \in (0,T) $ and every $\phi \in C^1([0,T] \times \overline{\Omega}; \R^d)$; while the energy inequality~\eqref{eq:mvenergyinequality} now reads

\begin{equation}\label{eq:mvenergyinequality2}
    \begin{aligned}
    \E_*^{mv}(\tau) \leq \E_*^{mv}(0) 
    &- \intO\intO\psi(x-y)\brk*{\skp*{\nu_{t,x};s|\vv|^2}\skp*{\nu_{t,y};s} - \skp*{\nu_{t,x};s\vv}\skp*{\nu_{t,y};s\vv}} \dx{x}\dx{y}
    \\[0.5em]
    &- \intO\intO \psi(x-y) \dx{m_\tau^{\vr|\vu|^2}}(x)\dx{m_\tau^\vr}(y)
    + \intO\intO \psi(x-y) \dx{m_\tau^\vr}(x)\dx{m_\tau^\vr}(y),
    \end{aligned}
\end{equation}
where
\begin{equation}
    \begin{aligned}
    \E_*^{mv}(\tau) &\coloneqq 
    \int_\Omega \skp*{\nu_{\tau,x}; \frac{1}{2} s | \vv|^2} \dx{x} + \frac12 m_\tau^{\rho|\vu|^2}(\Omega) 
    + \int_\Omega\skp*{ \nu_{\tau,x}; \frac12|\vf|^2} \dx{x} + \frac12m_\tau^{|\nabla\Phi|^2}(\Omega) 
    + \int_\Omega\skp*{ \nu_{\tau,x}; s} V(x) \dx{x}
	\\[0.5em]
	&\hspace{0.2cm}+\intO V(x)\dx{m}^{\rho}_\tau(x)
	+ \frac12\intO \skp*{\nu_{\tau,x};s}\intO W(x-y)\skp*{\nu_{\tau,y};s}\dx{y}\dx{x}
	+\frac12\intO\intO W(x-y) \dx{m_{\tau}^\vr}(y)\dx{m_{\tau}^\vr}(x)
	\\[0.5em]
	&= \intO \frac{1}{2} \ol{\rho |\vu|^2} + \frac{1}{2} \ol{|\nabla_x\Phi_\rho|^2} + \ol{\rho} V + \frac12\ol{\rho}\prt*{W\star\ol{\rho}} \dx{x}.
    \end{aligned}
\end{equation}

\bigskip

\noindent We wish to derive an analogue of Theorem~\ref{thm:mv-stronguniqness} for the system~\eqref{eq:euler_general}. To this end we again suppose that a $C^1$-regular solution $(r,\vU,\Phi_r)$ with $r>0$ is available, and consider a measure-valued solution with the same initial data $(r_0,\vU_0)$. Let us remark that the issue of existence of such a measure-valued solution can be settled by means of a similar two-step approximation argument as in Section~\ref{sec:existence}. Since the additional convolution-type terms pose no additional difficulties in providing existence of approximating sequences and passing to the limit, we skip the details.

Let us now discuss how the new terms influence the calculations towards a Gronwall inequality performed in the previous section. We keep the same relative energy functional as in~\eqref{eq:mvrelenergy}.
First, it can be readily seen that the $-\rho\nabla_x V$ term behaves in exactly the same way as the $-\rho\nabla\prt*{\frac12|x|^2}$ term previously. When we arrive at inequality~\eqref{eq:mvrelativeenergy}, we now have two additional lines on the right-hand side, namely
\begin{equation}
    \begin{aligned}\label{eq:W-line}
    -\frac12\brk*{\intO \ol{\vr}\prt*{W\star\ol{\vr}}\dx{x}}_{t=0}^{t=\tau} +\int_0^\tau\intO \ol{\vr}\prt*{\nabla_x W\star\ol{\vr}}\cdot\vU \dx{x}\dx{t} 
    + \int_0^\tau\intO \ol{\vr\prt*{\vu-\vU}}\cdot \nabla_x\prt*{W\star r}\dx{x}\dx{t},
    \end{aligned}
\end{equation}
and
\begin{equation}
    \begin{aligned}\label{eq:Psi-line}
    &-\int_0^\tau\intO\intO\psi(x-y)\ol{\vr}(t,x)\ol{\vr\vu}(t,y)\cdot\vU(t,x) \dx{y}\dx{x}\dx{t}
    +\int_0^\tau\intO\intO\psi(x-y)\ol{\vr}(t,y)\ol{\vr\vu}(t,x)\cdot\vU(t,x) \dx{y}\dx{x}\dx{t}
    \\[0.5em]
    &-\int_0^\tau\intO \ol{\vr\prt*{\vu-\vU}}\cdot\intO\psi(x-y)r(t,y)\prt*{\vU(t,y)-\vU(t,x)}\dx{y} \dx{x}\dx{t}
    \\[0.5em]
    &-\int_0^\tau\intO\intO\psi(x-y)\brk*{\ol{\vr|\vu|^2}(t,x)\ol{\vr}(t,y) - \ol{\vr\vu}(t,x)\ol{\vr\vu}(t,y)}\dx{y}\dx{x}\dx{t}.
    \end{aligned}
\end{equation}
For~\eqref{eq:W-line} we can, after some straightforward computations, write equivalently
\begin{equation}
    \begin{aligned}\label{eq:W-line2}
    -\frac12\brk*{\intO \ol{\vr}\prt*{W\star\ol{\vr}}\dx{x}}_{t=0}^{t=\tau} 
    + \int_0^\tau\intO\ol{\vr\vu}\cdot\nabla_x\prt*{W\star\ol{\vr}} \dx{x}\dx{t}
    + \int_0^\tau\intO\ol{\vr\prt*{\vu-\vU}}\cdot \nabla_x\prt*{W\star \prt*{r-\ol{\vr}}}\dx{x}\dx{t}.
    \end{aligned}
\end{equation}

Since the attraction-repulsion kernel $W$ is smooth, the first two terms of~\eqref{eq:W-line2} cancel according to the following calculation, see also~\cite[Section~5, Step~4.]{CaFeGwSw2015},
\begin{equation}
    \begin{aligned}
    -\frac12\brk*{\intO \ol{\vr}\prt*{W\star\ol{\vr}}\dx{x}}_{t=0}^{t=\tau}
    =
    -\frac12\int_0^\tau\intO \frac{\partial}{\partial t}\prt*{\ol{\vr}\prt*{W\star\ol{\vr}}}\dx{x}
    =
    -\int_0^\tau\intO\ol{\vr}\prt*{W\star\partial_t\ol{\vr}}\dx{x} = -\int_0^\tau\intO\ol{\vr\vu}\cdot\nabla_x\prt*{W\star\ol{\vr}}\dx{x}.
    \end{aligned}
\end{equation}

We now rewrite~\eqref{eq:Psi-line} as follows

\begin{equation}
    \begin{aligned}
    \int_0^\tau&\intO\intO \psi(x-y)\ol{\vr\prt*{\vu-\vU}}(t,x)\cdot\prt*{\vU(t,y)-\vU(t,x)}\ol{\prt*{\vr-r}}(t,y) \dx{y}\dx{x}\dx{t}
    \\[0.5em]
    &-\int_0^\tau\intO\intO \psi(x-y)\brk*{\ol{\vr\abs*{\vu-\vU}^2}(t,x)\ol{\vr}(t,y) - \ol{\vr\prt*{\vu-\vU}}(t,x)\cdot\ol{\vr\prt*{\vu-\vU}}(t,y)}\dx{y}\dx{x}\dx{t}.
    \end{aligned}
\end{equation}

The latter of the two terms can be discarded, since it is negative. Indeed, at the level of approximating sequences we can write, due to symmetry of $\psi$,
\begin{equation}
    \begin{aligned}
    &\intO\intO\psi(x-y)\brk*{\rho^\ep(t,y)\rho^\ep(t,x)\abs*{\vu^\ep(t,x)-\vU(t,x)}^2 - \rho^\ep(t,x)\prt*{\vu^\ep(t,x)-\vU(t,x)}\cdot\rho^\ep(t,y)\prt*{\vu^\ep(t,y)-\vU(t,y)}}\dx{y}\dx{x}\\[0.5em]
    &=\frac12\intO\intO\psi(x-y)\rho^\ep(t,x)\rho^\ep(t,y)\abs*{\prt*{\vu^\ep(t,x)-\vU(t,x)}-\prt*{\vu^\ep(t,y)-\vU(t,y)}}^2\dx{y}\dx{x}
    \geq 0.
    \end{aligned}
\end{equation}

\bigskip

\noindent It now only remains to bound the remaining two terms
\begin{equation}
    I_1 \equiv    \int_0^\tau\intO\ol{\vr\prt*{\vu-\vU}}\cdot \nabla_x\prt*{W\star \prt*{r-\ol{\vr}}}\dx{x}\dx{t},
\end{equation}
and
\begin{equation}
   I_2 \equiv \int_0^\tau\intO\intO \psi(x-y)\ol{\vr\prt*{\vu-\vU}}(t,x)\cdot\prt*{\vU(t,y)-\vU(t,x)}\ol{\prt*{\vr-r}}(t,y) \dx{y}\dx{x}\dx{t},
\end{equation}
in terms of the relative energy.
The strategy is to mimic the analogous bounds for the weak-strong case, as presented in~\cite{CaFeGwSw2015}. In fact, on the level of approximating sequences the calculations are exactly the same, since they only rely on functional inequalities and the Poisson equation. We present these arguments below for the readers' convenience.
First, notice that $I_1 = \lim_{\ep\to0}I_1^\ep$, where
\begin{equation}
    \begin{split}
        |I_1^\ep| &= \abs*{\int_0^\tau\intO\rho^\ep(t,x)\prt*{\vu^\ep(t,x)-\vU(t,x)}\cdot\intO\nabla_xW(x-y)\prt*{r(t,y)-\rho^\ep(t,y)}\dx{y}\,\dx{x}\dx{t}}
        \\[0.5em]
        &\leq c\int_0^\tau\intO\rho^\ep\abs*{\vu^\ep-\vU}^2\dx{x}\dx{t} + c\int_0^\tau\prt*{\intO\rho^\ep(t,x)\dx{x}}\norm{\nabla_xW\star\prt*{r(t,\cdot)-\rho^\ep(t,\cdot)}}^2_{L^\infty}\dx{t}
        \\[0.5em]
        &\leq c\int_0^\tau\intO\rho^\ep\abs*{\vu^\ep-\vU}^2\dx{x}\dx{t} + cM_{\rho^\ep}\int_0^\tau\norm{r(t,\cdot)-\rho^\ep(t,\cdot)}^2_{W^{-1,2}}\dx{t}
        \\[0.5em]
        &\leq c\int_0^\tau\intO\rho^\ep\abs*{\vu^\ep-\vU}^2\dx{x}\dx{t} + cM_{\rho^\ep}\int_0^\tau\norm{\nabla_x\Phi_r(t,\cdot)-\nabla\Phi_{\rho^\ep}(t,\cdot)}^2_{L^2}\dx{t}.
    \end{split}
\end{equation}
In the limit $\ep\to0$ we obtain the desired estimate
\begin{equation}
    I_1 \leq C\int_0^\tau\E_{rel}^{mv}(t)\dx{t}.
\end{equation}
The second term, $I_2$, is treated identically once we use boundedness of the strong solution $\vU$. 

\bigskip

We can therefore once again infer that whenever the regular and the measure-valued solution emanate from the same initial data, their relative energy vanishes for a.a.\ times $t>0$. 
As in the previous section we deduce that the projection of the Young measure $\nu_{t,x}$ onto the third coordinate is the Dirac measure concentrated at $\nabla_x\Phi_r(t,x)$, and $m^{|\nabla\Phi|^2}=0$. 
Using this information in the Poisson equation, we obtain that $\ol{\rho}=r$ almost everywhere. Similarly, from the vanishing of the first term in the relative energy,
\begin{equation}
    \intO\ol{\rho\abs*{\vu-\vU}}\dx{x} = 0,
\end{equation}
we deduce the identifications
\begin{equation}
	\ol{\rho|\vu|^2} = r|\vU|^2,\quad \ol{\rho\vu\otimes\vu} = r\vU\otimes\vU.
\end{equation}
Substituting these identifications into the (measure-valued) momentum equation of~\eqref{eq:euler_general} and using the strong formulation for $(r,\vU)$, we obtain the following ODE
\begin{equation}
    \partial_t\prt*{\ol{\rho\vu} - r\vU} = r\brk*{\psi\star\prt*{\ol{\rho\vu}-r\vU}} - \prt*{\psi\star r}\prt*{\ol{\rho\vu}-r\vU},
\end{equation}
which readily implies that $\ol{\rho\vu}=r\vU$ almost everywhere in $(0,T)\times\Omega$. Whence we obtain the following result

\begin{theorem}
    Let $1\leq d\leq 3$ and $\Omega\subset\R^d$ be a bounded smooth domain. Let 
    \[
        (r,\vU,\Phi_r)\in C^1([0,T)\times\bar{\Omega};(0,\infty))\times C^1([0,T)\times\bar{\Omega};\R^d)\times C^2([0,T)\times\bar{\Omega})
    \]
    be a strong solution of~\eqref{eq:euler_general} with initial data $r(0,x)=r_0(x),\, \vU(0,x)=\vU_0(x)$ of finite energy, and let\\ $(\vnu, m^\rho, m^{\rho\vu}, m^{\rho \vu\otimes \vu}, m^{|\nabla\Phi|^2}, m^{\nabla\Phi\otimes\nabla\Phi})$ be a dissipative measure-valued solution with initial state
    \[
        \nu_{0,x} = \delta_{\{r_0,\vU_0,\nabla\Phi_r(0,x)\}}\;\;\;\; \text{for a.e.}\;\; x\in\Omega.
    \]
Then 
\[
m^{\nabla\Phi\otimes\nabla\Phi}=0,\;\; m^{|\nabla\Phi|^2}=0,
\]
and we have the following identifications
\begin{align}
	\skp*{\nu_{t,x} ; \rho} + m^{\rho} &= r,\\	
	\skp*{\nu_{t,x} ; \rho\vu} + m^{\rho\vu} &= r\vU,\\
	\skp*{\nu_{t,x} ; \rho\vu\otimes\vu} + m^{\rho\vu\otimes\vu} &= r\vU\otimes\vU,\\
	\skp*{\nu_{t,x} ; \rho|\vu|^2} + m^{\rho|\vu|^2} &= r|\vU|^2,
\end{align}
which hold for almost every $(t,x)\in(0,T)\times\R^d$.
Furthermore, the Young measure admits the decomposition
\begin{equation}
	\nu_{t,x} = \bar{\nu}_{t,x} \otimes \delta_{\{\nabla\Phi_r(t,x)\}},
\end{equation}
for some parameterised measure $\bar{\nu}\in L^\infty_{\weak}((0,T)\times\R^d;\M^+([0,\infty)\times\R^d))$; and in turn the restriction $\bar{\nu}\mres ((0,\infty)\times\R^d)$ decomposes into
\begin{equation}
	\bar{\nu}_{t,x} \mres ((0,\infty)\times\R^d) = \bar{\bar{\nu}}_{t,x} \otimes \delta_{\{\vU(t,x)\}}
\end{equation}
for some parameterised measure $\bar{\bar{\nu}}\in L^\infty_{\weak}((0,T)\times\R^d;\M^+(0,\infty))$.
Finally, all the non-zero concentration measures $m^\rho$, $m^{\rho\vu}$, $m^{\rho\vu\otimes\vu}$, $m^{\rho|\vu|^2}$ are absolutely continuous with respect to the Lebesgue measure.
\end{theorem}


\begin{appendices}

\section{Appendix}
\label{sec:appendix}

\subsection{Young measures}
\label{sec:appA}

Below we gather some additional facts about the parameterised measure generated by our approximating sequences of solutions which we used to pass to the limit in Section~\ref{sec:existence} and deduce the weak-strong identifications in Section~\ref{sec:mainproof}. The required notation and definitions are presented in the introduction.

\begin{lemma}
\label{lem:integrability}
	Suppose $X\subset\R^n$ is bounded. Let $z^\ep:X\to Y$ be a sequence of measurable functions and let $\vnu=(\nu_x)\in L^\infty_{\weak}(X;\M^+(Y))$ denote the assosciated Young measure. Let $f\in C(Y)$ be a continuous function and suppose that the sequence $(f(z^\ep))$ is uniformly bounded in $L^1(X)$, \ie,
	\begin{equation}
    	\sup_{\ep>0}\int_X |f(z^\ep(x))|\dx{x} \leq C.
	\end{equation}
Then the function $f$ is $\vnu$-measurable, \ie, the map $x\mapsto\skp*{\nu_x ; f}$ is well-defined for a.e.\ $x\in X$. Moreover, the map $x\mapsto\skp*{\nu_x ; f}$ belongs to $L^1(X)$. 	
\end{lemma}

\begin{proof}
	Without loss of generality we can assume that $f\geq 0$.
	Integration with respect to the Young measure is well-defined for continuous functions which vanish at infinity. So consider the sequence of truncated functions $f^k(y) = \theta^k(|y|) f(y)$ where
	\begin{equation}
		\theta^k(\alpha)=
			\begin{cases}
				1 & \text{if $|\alpha| < k$}\\
				(k+1)-\alpha & \text{if $k\leq |\alpha| \leq k+1$}\\
				0 & \text{if $|\alpha|>k+1$}.
			\end{cases}
	\end{equation}
Then $f^k\in C_0(Y)$, $0\leq f^k\leq f$, the seqeunce is non-decreasing and converges pointwise to $f$. It follows from the Monotone Convergence Theorem that $\skp*{\nu_x; f}$ is well-defined for a.e.\ $x\in X$. Moreover, for each $\phi\in L^1(X)$ we have
	\begin{equation}
	\int_X \phi(x)f^k(z^\ep(x))\dx{x} \longrightarrow \int_X \skp*{\nu_x;f^k}\phi(x)\dx{x}.
	\end{equation}
In particular, for $\phi\equiv 1$ we have
	\begin{equation}
	\int_X f^k(z^\ep(x))\dx{x} \longrightarrow \int_X \skp*{\nu_x;f^k} \dx{x}.
	\end{equation}
But
	\begin{equation}
		\int_X f^k(z^\ep(x))\dx{x} \leq \int_X f(z^\ep(x))\dx{x} \leq C,
	\end{equation}
and so the integrals
	\begin{equation}
		\int_X \skp*{\nu_x;f^k} \dx{x} \leq C
	\end{equation}
are bounded uniformly in $k$.
Therefore, by monotone convergence, we deduce that $\brk*{x\mapsto\skp*{\nu_x;f}} \in L^1(X)$.	
\end{proof}

Consequently, defining the concentration-defect measure as in~\eqref{eq:concentrationdefect}, we can deduce the convergence
\begin{equation}
\label{eq:limitwithconcentration}
	\int_X f(z^\ep(x))\phi(x)\dx{x} \longrightarrow \int_X \skp*{\nu_x ; f}\phi(x)\dx{x} + \int_X \phi(x)\dx{m^f}(x)
\end{equation}
for every bounded continuous test function $\phi\in C({\ol{X}})$. The reader might like to compare the above representation result with the notion of biting convergence and how the Young measure describes the biting limit~\cite{BallMurat1989, Pedregal}.

\bigskip

The next result concerns the possible problem with canonical projections of the parameterised measure onto one of the components of the dummy vector $(s,\vv,\vf)$. This another technical issue stemming from lack of tightness. Namely, the projection of the Young measure generated by a multi-component sequence onto one of its coordinates might not agree with the Young measure generated by the coresponding component. For example, consider the sequence $z^n(x) = (1,n)$ in $\R^2$. Then the corresponding Young measure $\nu$ is zero almost everywhere, while the Young measure generated by the projection sequence $\pi_1z^n =1$ is equal to $\delta_{\{1\}}$. This effect cannot occur if the Young measure generated by $z^n = (z_1^n, z_2^n)$ is a probability measure. Indeed, suppose this is the case. Then taking $f(z^n) = f_1(z_1^n)$ for $f_1\in C_0$, we have
\begin{equation}
	\int f_1(z_1^n)\phi = \int f(z^n)\phi \longrightarrow \int \skp*{\nu_x ; f} \phi = \int \skp*{\pi_1\nu_x ; f_1}\phi,
\end{equation}
and
\begin{equation}
	\int f_1(z_1^n)\phi \longrightarrow \int \skp*{\eta_x;f_1}\phi,
\end{equation}
where $\eta$ denotes the Young measure generated by the sequence $(z_1^n)$. It follows that $\nu_x=\eta_x$ almost everywhere.
\\
In the current context the above issue could potentially lead to problems when considering the Young measure $\vnu$ coming from an approximating sequence $z^\ep=(\rho^\ep, \vu^\ep)$. This is remedied by the observation that on each set $\{\rho^\ep\geq\alpha\}$, $\alpha>0$, the sequence $z^\ep$ satisfies the tightness condition, while on the vacuum zones we are free to modify the measures in question arbitrarily. More precisely, we have
\begin{lemma}
\label{lem:projections}
	Let $z^\ep = (\rho^\ep, \vu^\ep)$ be any sequence of approximate solutions such that $\rho^\ep$ and $\rho^\ep\vu^\ep$ are uniformy bounded in $L^1(\R^{d+1})$. Let $\vnu$ be the Young meaure generated by $z^\ep$ and $\vec{\eta}$ be the Young measure generated by $\rho^\ep$. Then
	\begin{equation}
		\pi_1\prt*{\nu_x\mres\prt*{(0,\infty)\times\R^d}} = \eta_x\mres(0,\infty)
	\end{equation}
for almost every $x\in\R^{d+1}$.	
\end{lemma}

\begin{proof}
	Let $\alpha>0$ be fixed. Let $\vec{\mu}_\alpha$ denote the Young measure generated by $z^\ep$ considered on the set $\{\rho^\ep>\alpha\}$. Then $\vec{\mu}\in L^\infty_{\weak}(\R^{d+1};\PP((\alpha,\infty)\times\R^d))$. Clearly $\vec{\mu}_\alpha = \vec{\mu}\mres\prt*{(\alpha,\infty)\times\R^d}$. Similarly define $\vec{\eta}_\alpha = \vec{\eta}\mres(\alpha,\infty)$.
Then, as discussed above, we have $\pi_1\vec{\mu}_\alpha = \vec{\eta}_\alpha$.\\
Now choose a Borel set $A\subset(0,\infty)$ and denote $A^\alpha = A\cap(0,\alpha)$. Then, by definition, $\vec{\eta}_\alpha(A) = \vec{\eta}(A^\alpha)$. Since the family $A^\alpha$ is non-decreasing as $\alpha\to0$ and $A = \bigcup_{\alpha\geq 0}A^\alpha$, we have $\vec{\eta}(A)=\lim_{\alpha\to0}\vec{\eta}(A^\alpha)$. Therefore the sequence $(\vec{\eta}_\alpha(A))_{\alpha>0}$ of extended reals converges as $\alpha\to0$ to $\eta(A)$ for every Borel set $A$.\\
Similarly, $\pi_1\vec{\mu}_\alpha(A) = \vec{\mu}(A\times\R^d) = \lim_{\alpha\to0}\vec{\mu}(A^\alpha\times\R^d) = \lim_{\alpha\to0}\pi_1\vec{\mu}(A^\alpha)$.
But on the other hand
$\vec{\mu}(A^\alpha\times\R^d) = \vec{\mu}_\alpha(A\times\R^d) = \pi_1\vec{\mu}_\alpha(A)$, and thus the sequence $(\pi_1\vec{\mu}_\alpha(A))_{\alpha>0}$ converges to $\pi_1\vec{\mu}(A)$.
Consequently, we must have $\vec{\eta}(A) = \pi_1\vec{\mu}(A)$ for every Borel set $A\subset(0,\infty)$.
\end{proof}

\subsection{Inequalities}
\label{sec:appB}
We provide here the statement and proof of the simple geometric inequalities which we used in the proof of the main theorem to identifiy certain weak limits.

\begin{lemma}
\label{lem:inegalites}
	Let $\vU\in L^\infty(\R^{n};\R^d)$ be a bounded vector-valued function. Then for any vector $\vu\in\R^d$ and $\delta>0$ small enough we have the following inequalities
	\begin{enumerate}[align=left]
		\item 
				$\displaystyle\abs*{u_iu_j - U_iU_j} \leq c\delta + C_\delta\prt*{|u_i-U_j|^2+|u_j-U_j|^2},\quad\text{for any }\; i,j=1,\dots,d$,
		\item
				$\displaystyle \abs*{|\vu|^2 - |\vU|^2} \leq c\delta + C_\delta|\vu-\vU|^2$
	\end{enumerate}
	where the positive constants $c$ and $C_\delta$ depend only on $\norm{\vU}_\infty$ and $\delta$ and $\norm{\vU}_\infty$, respectively.
\end{lemma}

\begin{proof}
Firstly, we fix a point $y\in\R^{n}$ and consider the fixed vector $\vU = \vU(y)$; furthermore we choose $0<\delta<1$.
Let us denote $p(\vu) = u_iu_j-U_iU_j$ and consider the change of variables $v_i = u_i-U_i$ and $v_j=u_j-U_j$. Then we have
\begin{equation}
	p(\vv) = v_iv_j + U_jv_i + U_iv_j.
\end{equation}
Now we observe that whenever $\min(|v_i|,|v_j|)>1$ we have
\begin{equation}
	|p(\vv)| \leq (v_i^2 + v_j^2) + \norm{\vU}_{\infty}v_i^2 + \norm{\vU}_{\infty}v_j^2 \leq C(v_i^2 + v_j^2);
\end{equation}
while whenever $\max(|v_i|,|v_j|)>1$, then
\begin{equation}
	|p(\vv)| \leq (v_i^2 + v_j^2) + 2\norm{\vU}_{\infty}\max(v_i^2, v_j^2) \leq C(v_i^2 + v_j^2);
\end{equation}
and whenever $\max(|v_i|,|v_j|)\leq\delta$, then
\begin{equation}
	|p(\vv)| \leq \delta^2 + 2\norm{\vU}_{\infty}\delta \leq c\delta.
\end{equation}
In the remaining cases, we use continuity of the polynomial $p$:\ for the compact sets
\begin{equation}
	X_1 = \{\delta\leq |v_i|, |v_j| \leq 1\},\;\;\; X_2 = \{\delta\leq |v_i| \leq1 , |v_j| \leq \delta\},\;\;\; X_3 \{ |v_i|\leq \delta, \delta \leq |v_j| \leq 1\}
\end{equation}
there are finite constants (for $\alpha = 1,2,3$)
\begin{equation}
	K^\alpha_\delta = \sup_{X_\alpha}\;|p(\vv)| < \infty,
\end{equation}
so that
\begin{equation}
	|p(\vv)| \leq K^\alpha_\delta \leq \frac{K^\alpha_\delta}{\delta^2}\max(v_i^2, v_j^2) \leq C_\delta(v_i^2 + v_j^2).
\end{equation}
Altogether we obtain the first of the claimed inequalities.
The second one follows from taking $i=j$ and summing over all indices $i=1,\dots,d$.
\end{proof}

\begin{lemma}
\label{lem:inegalites2}
	Let $\vv\in\R^d$ and consider the matrix $A = \vv\otimes \vv$. Then the $L^1$-norm of $A$ can be bounded by its trace:
	\begin{equation}
		|A| = \sum_{i,j=1}^d |a_{ij}|  \leq c_d\; \mathrm{tr}(A) = c_d |\vv|^2.
	\end{equation}
\end{lemma}
\begin{proof}
	From the elementary inequality
	\begin{equation}
		|v_iv_j| \leq \frac12\prt*{|v_i|^2 + |v_j|^2},
	\end{equation}
	we have the bound
	\begin{equation}
		|a_{ij}| \leq \frac12(|a_{ii}| + |a_{jj}|).
	\end{equation}
\end{proof}
This elementary inequality is used for instance in the calculation towards the relative energy inequality to bound the oscillation and concentration parts of the quadratic term $\rho(\vu-\vU)\otimes(\vu-\vU)$ by the "kinetic" term $\rho|\vu-\vU|^2$ of the relative energy. Also, it is used to deduce that $|m^{\vv\otimes \vv}| \leq m^{|\vv|^2}$.
\end{appendices}

\section*{Acknowledgments}

JAC was partially supported by EPSRC grant number EP/P031587/1 and the Advanced Grant Nonlocal-CPD (Nonlocal PDEs for Complex Particle Dynamics: Phase Transitions, Patterns and Synchronization) of the European Research Council Executive Agency (ERC) under the European Union's Horizon 2020 research and innovation programme (grant agreement No. 883363).
TD was partially supported by National Science Centre, Poland, under agreement no UMO-2018/31/N/ST1/02394, the Polish National Agency of Academic Exchange (NAWA), and the Foundation for Polish Science. A\'{S}-G and PG were supported by National Science Centre, Poland, under agreement no UMO-2017/27/B/ST1/01569.

\end{document}